\documentstyle[11pt]{article}

\begin{document}
\noindent
\begin{center}
  {\LARGE Stringy geometry and topology of orbifolds}\footnote{ partially
supported by the National Science Foundation}
  \end{center}

  \noindent
  \begin{center}

     {\large  Yongbin Ruan}\\[5pt]
      Department of Mathematics, University of Wisconsin-Madison\\
        Madison, WI 53706\\[5pt]
        Email:  ruan@math.wisc.edu\\[5pt]

        Dedicated to the occasion of Herb Clemens' 60-th
        birthday\\
              \end{center}

              \def \x{{\bf x}}
              \def \J{{\cal J}}
              \def \M{{\cal M}}
              \def \A{{\cal A}}
              \def \B{{\cal B}}
              \def \C{{\bf C}}
              \def \Z{{\bf Z}}
              \def \R{{\bf R}}
              \def \P{{\bf P}}
              \def \I{{\bf I}}
              \def \N{{\cal N}}
              \def \T{{\cal T}}
              \def \Q{{\bf Q}}
              \def \D{{\cal D}}
              \def \H{{\cal H}}
              \def \S{{\cal S}}
              \def \e{{\bf E}}
              \def \CP{{\bf CP}}
              \def \U{{\cal U}}
              \def \E{{\cal E}}
              \def \F{{\bf C}}
              \def \L{{\cal L}}
              \def \K{{\cal K}}
              \def \G{{\cal G}}
              \def \z{{\bf z}}
              \def \m{{\bf m}}
              \def \n{{\bf n}}
              \def \V{{\cal V}}
              \def \W{{\cal W}}
 \def \g{{\bf g}}
              \def \h{{\bf h}}

\tableofcontents

\section{Introduction}
In 1985, Dixon, Harvey, Vafa and Witten considered  string
    theory on orbifolds (arising as global
quotients $X/G$ by a finite group $G$)\cite{DHVW}. Although an
orbifold is a singular space, orbifold string theory is
surprisingly  a smooth theory. Since then, orbifold string theory
has become a rather important part of the landscape of string
theory. A search in hep-th yields more than 200 papers whose title
contains the word "orbifold".  Although orbifold string theory has
been around for a while, apparently it was poorly explored in
geometry and topology. For last fifteen years, only a small piece
of orbifold string theory concerning orbifold Euler-Hodge numbers
has been studied in geometry and topology in relation to the McKay
correspondence in algebraic geometry (see \cite{Re}, \cite{B1}).
The highlight of previous works was Batyrev's proof that orbifold
Hodge numbers of a global quotient $X/G$ for a finite subgroup
$G\subset SL(n, \C)$ are the same as Hodge numbers of its crepant
resolution. We shall see later that this is only a special case of
many "Orbifold string theory conjectures".

    One of the reasons for such a slow development of the mathematical aspect of
    orbifold string theory is  poor communication between mathematicians and physicists.
    The situation has been vastly improved in the last few years due to the
    increasing number of mathematicians who have been  making a serious effort to
    understand physics. With our improved understanding of the physical ideas
    behind orbifold string theory, we are starting to be able to explore the full
    implications of orbifold string theory. In the last year, a series of work
    \cite{CR1,CR2}, \cite{R}, \cite{AR} have been carried out on this direction.

    Even with a superficial understanding of orbifold string
    theory, it is obvious that the mathematics surrounding
    orbifold string theory must be striking.
    In fact, it has motivated so much new mathematics unique to
    orbifolds. We believe that there is an emerging  "stringy" topology and geometry of orbifolds.
    The core of this new geometry and topology is the concept of
    twisted sectors. Roughly speaking, the consistency of orbifold string theory requires that
     the string
Hilbert space has to contain factors called twisted sectors.
Twisted sectors can be viewed as the contribution from
singularities. All other quantities such as correlation functions
have to contain the contributions from the twisted sectors. In
    another words, the ordinary topology of orbifolds is a WRONG
    theory. The correct one must incorporate twisted sectors. Furthermore, orbifold string theory
    has a certain internal  freedom (discrete torsion) \cite{V}. Discrete torsion will allow
    us to twiste orbifold string theory \cite{VW}. These are the most important
    new conceptual ingredients in orbifold string theory. We will
    emphasis them in our mathematical construction as well.

    Once   we overcome these basic conceptual hurdles, we can explore the implications of
    the inclusion of twisted sectors and discrete torsion in all aspect of topology and geometry. This is what I would like
    call "Stringy geometry and topology of orbifold". Right now, this new subject is very much in
    its infancy. Currently, most of our motivation comes from physics. As time goes on,
    I expect that more mathematical motivations will emerge.

    In this article, we will survey the new developments on this
    subject. For a new subject, it is common that there are more
    problems and speculations than the mathematics we can actually prove.
    This     is also the case for the present subject. Therefore, we will also
    spend considerable time in talking about problems and
    conjectures.

    One of the major problems of this subject is the lack of references. In section 2,
    we will give a self-contained introduction to
    orbifolds. Furthermore, we will introduce the key technical
    concept of good map.

    In section 3, I will carry out the construction of orbifold cohomology
    \cite{CR1}, \cite{R}. Orbifold string theory has an internal
     freedom (discrete torsion), which will allow us to twist
    orbifold cohomology. However, discrete torsion is not enough account for
    all the known examples. In \cite{R}, a more general notion of inner
    local system was introduced. We shall construct our orbifold
    cohomology in this general setting.

    In section 4, we will change our point of view to K-theory
    and develop orbifold K-theory. Again, we will incorporate
    discrete torsion in our theory. Here, the mathematical
    motivation is projective representations. We will give a
    detailed description of this approach and establish the additive
    isomorphism between orbifold K-theory and orbifold cohomology.
    A surprising byproduct of orbifold K-theory is the new
    multiplicative structure between DIFFERENT twisted orbifold
    cohomologies, which is impossible to observe from the cohomological
    point of view. Right now, this new multiplicative structure is
    very intriguing. It is certainly worth more investigation.

    In section 5, we will shift from classical theory to
    quantum theory---orbifold quantum cohomology \cite{CR2}. We will
    introduce  the notion of orbifold stable map, which is a
    nontrivial generalization of stable maps. Then, we will study
    various properties of orbifold stable maps and construct
    orbifold quantum cohomology.

    In section 6, we will focus on some of the main predictions
    from orbifold string theory. The main idea is that orbifold
    theory should predict the ordinary theory of its
    desingularizations. A desingularization $Y$ of a Gorenstein
    orbifold $X$ is obtained by deforming $X$ and then taking a
    crepant resolution. Two extreme cases are the ones obtained by
    either deformation or resolution alone. There is a body of
    conjectures about their relations depending on the particular
    setting. We call all of  them by the term "K or Q-orbifold string theory
    conjecture". Here, the letter K or Q is indicating the particular setting
    we are talking about. Furthermore, it can be combined further with
    author's quantum minimal model conjecture to extend the orbifold
    string theory conjecture and the quantum minimal model conjecture
    in a natural way. Other problems will also be discussed.

    Finally, a historical note is in order. So far, the physical
    construction of orbifold string theory has only been carried out on
    global quotients of the form $X/G$. It is not clear how to
    do it over general orbifolds. However, most of the important
    examples such as Calabi-Yau orbifolds are not global
    quotients. In
    dimension three every Calabi-Yau orbifold admits a crepant
    resolution. Hence, we can restrict ourselves to smooth models. In higher dimension, this is no longer
    true.
     If
    we want to extend  wonderful theories such as mirror
    symmetry to higher dimension, we are forced to work over
    singular manifolds. Therefore, it is necessary to be able to construct
    orbifold string theory over general orbifolds. It is our hope that a better mathematical
    understanding of general orbifolds will help the physical
    construction as well. Throughout the paper, we will put an
    emphasis on developing the theory over general orbifolds.

    Even
    in the case of global quotients, the best understood part of orbifold string theory
    is orbifold conformal field theory. It is not clear how to do
    new geometry and topology except the orbifold Hodge number \cite{Z}.
    However, there are many papers about McKay correspondences for
    global quotients. One can find relevant references in \cite{Re}.
    An equivalent formulation of orbifold stable maps was studied in algebraic geometry
    independently by D. Abramovich and A. Vistoli \cite{AV}.

    I would like to express my special thanks to my collaborators
    Alejandro Adem and Weimin Chen. I also would like to thank R.
    Dijkgraaf, E. Witten and E. Zaslow  for stimulating
    discussions about orbifold string theory.

\section{Basics}
    One of the difficult of this subject is that there are few
    references. A detailed description of basic material  on orbifold
    has been given in the appendix of \cite{CR3}, which includes the
    crucial new concept of good map. Here we review the basic
    construction. Furthermore, we take a slightly more general definition
    of orbifold, which corresponds to a smooth Deligne-Mumford stack.
    It seems to be a convenient and  natural category even if we work with a more
    restrictive definition, which we call reduced orbifold.

\subsection{Orbifolds}
     Primary examples of orbifold  are quotient space of smooth manifolds by a
smooth finite group action. Here we imagine that the quotient
space is uniformized (or modeled) by the manifold with the finite
group action. We do not require  the group action to be effective.
If it happens to be effective, we call it a {\em reduced
orbifold}. It is clear that we can canonically associate a reduced
orbifold to an orbifold by redefining the group. Furthermore, if
an element acts nontrivially, we require that the fixed-point set
is of codimension at least two. This is the case, for example,
when the action is orientation-preserving. This requirement has
the consequence that the non-fixed-point set is locally connected.

\vspace{2mm}

Let $U$ be a connected topological space, $V$ be a connected
$n$-dimensional smooth manifold and $G$ be a finite group acting
on $V$ smoothly. {\it An n-dimensional uniformizing system} of $U$
is a triple $(V,G,\pi)$, where $\pi: V \rightarrow U$ is a
continuous map inducing a homeomorphism between $V/G$ and $U$. Two
uniformizing systems $(V_i,G_i,\pi_i)$, $i=1,2$, are {\it
isomorphic} if there is a diffeomorphism $\phi: V_1\rightarrow
V_2$ and an isomorphism $\lambda: G_1\rightarrow G_2$ such that
$\phi$ is $\lambda$-equivariant, and $\pi_2\circ\phi=\pi_1$. It is
easily seen that if $(\phi,\lambda)$ is an automorphism of
$(V,G,\pi)$, then there is a  $g\in G$ such that $\phi(x)=g\cdot
x$ and $\lambda(a)=g\cdot a\cdot g^{-1}$ for any $x\in V$ and
$a\in G$.

Let $i:U^\prime\hookrightarrow U$ be a connected open subset of
$U$, and $(V^\prime,G^\prime, \pi^\prime)$ be a uniformizing
system of $U^\prime$. We say that $(V^\prime,G^\prime,\pi^\prime)$
is induced from  $(V,G,\pi)$ if there is a monomorphism $\lambda:
G^\prime\rightarrow G$ and a $\lambda$-equivariant open embedding
$\phi: V^\prime\rightarrow V$ such that $i
\circ\pi^\prime=\pi\circ\phi$. We follow Satake [S] and call
$(\phi,\lambda): (V^\prime,G^\prime,\pi^\prime)\rightarrow
(V,G,\pi)$ an {\it injection}. Two injections $(\phi_i,\lambda_i):
(V^\prime_i,G^\prime_i,\pi^\prime_i)\rightarrow (V,G,\pi)$,
$i=1,2$, are {\it isomorphic} if there is an isomorphism
$(\psi,\tau)$ between $(V^\prime_1,G^\prime_1,\pi^\prime_1)$ and
$(V^\prime_2,G^\prime_2,\pi^\prime_2)$, and an automorphism
$(\bar{\psi}, \bar{\tau})$ of $(V,G,\pi)$ such that
$(\bar{\psi},\bar{\tau})\circ
(\phi_1,\lambda_1)=(\phi_2,\lambda_2)\circ (\psi,\tau)$.

\vspace{2mm}

\noindent{\bf Lemma 2.1.1: }{\it Let $(V,G,\pi)$ be a uniformizing
system of $U$. For any connected open subset $U^\prime$ of $U$,
$(V,G,\pi)$ induces a unique isomorphism class of uniformizing
systems of $U^\prime$.}

\vspace{2mm}

\noindent{\bf Proof:}

{\it Existence}: Consider the preimage $\pi^{-1}(U^\prime)$ in
$V$. $G$ acts as permutations on the set of connected components
of $\pi^{-1}(U^\prime)$. Let $V^\prime$ be one of the connected
components of $\pi^{-1}(U^\prime)$, $G^\prime$ be the subgroup of
$G$ which fixes the component $V^\prime$ and
$\pi^\prime=\pi|_{V^\prime}$. Then
$(V^\prime,G^\prime,\pi^\prime)$ is an induced uniformizing system
of $U^\prime$.

{\it Uniqueness}: First of all, different choices of the connected
components of $\pi^{-1}(U^\prime)$ induce isomorphic uniformizing
systems. Secondly, let $(V^\prime_1,G^\prime_1,\pi^\prime_1)$ be
any induced uniformizing system of $U^\prime$ and $(\psi,\tau)$ be
the injection into $(V,G,\pi)$. We will show that $(\psi,\tau)$
induces an isomorphism between
$(V^\prime_1,G^\prime_1,\pi^\prime_1)$ and an induced uniformizing
system given by a connected component of $\pi^{-1}(U^\prime)$.
Suppose $\psi( V^\prime_1)$ lies in the connected component
$V^\prime$. We can show that $\psi(V^\prime_1)$ is closed in
$V^\prime$. Let $\psi(x_n)\rightarrow y_0$ in $V^\prime$, $x_n\in
V^\prime_1$, then there exists a $z_0\in V^\prime_1$ such that
$\pi^\prime_1(z_0)=\pi(y_0)$, and $z_n\in V^\prime_1$ such that
$z_n\rightarrow z_0$,
 $\pi^\prime_1(z_n)=\pi(\psi(x_n))=\pi^\prime_1(x_n)$. So  there exist
$a_n\in G^\prime_1$ such that $a_n(z_n)= x_n $. Since $G^\prime_1$
is finite, it follows that for large $n$, $a_n=a$ is a constant.
So $x_n\rightarrow a(z_0)$ in $V^\prime_1$ and $y_0=\psi(a(z_0))$,
i.e., $\psi(V^\prime_1)$ is closed in $V^\prime$.
$\psi(V^\prime_1)$ is also open in $V^\prime$. So $\psi$ induces a
diffeomorphism between $V^\prime_1$ and $V^\prime $. From this we
can easily see that $(V^\prime_1,G^\prime_1,\pi^\prime_1)$ and
$(V^\prime, G^\prime,\pi^\prime)$ are isomorphic. \hfill $\Box$

\vspace{2mm}

Let $U$ be a connected and locally connected topological space.
For any point $p\in U$, we can define the {\it germ} of
uniformizing systems at $p$ in the following sense. Let $(V_1,
G_1,\pi_1)$ and $(V_2,G_2,\pi_2)$ be uniformizing systems of
neighborhoods $U_1$ and $U_2$ of $p$. We say that $(V_1,
G_1,\pi_1)$ and $(V_2,G_2,\pi_2)$
 are {\it equivalent} at $p$ if they induce
isomorphic uniformizing systems for a neighborhood $U_3$ of $p$.

\vspace{2mm}

\noindent{\bf Definition 2.1.2: }{\it Let $X$ be a Hausdorff,
second countable topological space. An {\it n-dimensional orbifold
structure} on $X$ is given by the following data: for any point
$p\in X$, there is a neighborhood $U_p$ and an $n$-dimensional
uniformizing system $(V_p, G_p, \pi_p)$ of $U_p$ such that for any
point $q\in U_p$, $(V_p, G_p, \pi_p)$ and $(V_q,G_q,\pi_q)$ are
equivalent at $q$ (i.e., define the same germ at $q$). The {\it
germ} of orbifold  structures on $X$ is defined in the following
sense: two orbifold structures $\{(V_p,G_p,\pi_p): p\in X\}$ and
$\{(V^\prime_p,G^\prime_p,\pi^\prime_p): p\in X\}$ are {\it
equivalent} if for any $p\in X$, $(V_p,G_p,\pi_p)$ and
$(V^\prime_p,G^\prime_p,\pi^\prime_p)$ are equivalent at $p$. With
a given germ of orbifold structures on it, $X$ is called an {\it
orbifold}. We call each $U_p$ a {\it uniformized neighborhood } of
$p$, and $(V_p,G_p,\pi_p)$ a {\it chart} at $p$. An open subset
$U$ of $X$ is called a {\it uniformized open set} if it is
uniformized by $(V,G,\pi)$ such that for each $p\in U$,
$(V,G,\pi)$ defines the same germ as $(V_p,G_p,\pi_p)$ at $p$.
 A point
$p\in X$ is called {\it regular or smooth} if $G_p$ is trivial;
otherwise, it is called {\it singular}. The set of smooth points
is denoted by $X_{reg}$, and the set of singular points is denoted
by $\Sigma X$. An orbifold $X$ is called a reduced orbifold if
$G_p$ acts effectively on $V_p$. In this case, we can choose $G_p$
to be a subgroup of $O(n)$.}

\vspace{2mm}
    \noindent
    {\bf Remark 2.1.3: }{\it It is obvious that we can associate
    canonically a reduced orbifold structure $X_{red}$ to each orbifold $X$ by redefining the local
    group $G'_p$ to be the quotient of $G_p$ divided by the
    subgroup fixing $V_p$.}
    \vskip 0.1in

\noindent{\bf Remark 2.1.4: }{\it There is a notion of {\it
orbifold with boundary}, in which we allow the uniformizing
systems to be smooth manifolds with boundary, with a finite group
action preserving the boundary. If $X$ is an orbifold with
boundary, then it is easily seen that the boundary $\partial X$
inherits an orbifold structure from $X$ and becomes an orbifold.}
\vskip 0.1in

\noindent{\bf Example 2.1.5: }{\it Let's consider the
2-dimensional sphere $S^2$. Let $D_s$, $D_n$ be open disc
neighborhoods of the south pole and the north pole such that
$S^2=D_s\cup D_n$. Let $D_s$ be uniformized by
$(\tilde{D}_s,\Z_2,\pi_s)$, and $D_n$ be uniformized by
$(\tilde{D}_n,\Z_3,\pi_n)$ where $\Z_2$, $\Z_3$ act on
$\tilde{D}_s$ and $\tilde{D}_n$ by rotations. For any point in
$S^2$ other than the south pole and the north pole, we take a
chart at it induced by either $(\tilde{D}_s,\Z_2,\pi_s)$ or
$(\tilde{D}_n,\Z_3,\pi_n)$. It is easily seen that this defines a
2-dimensional orbifold structure on $S^2$. Note that as an open
subset of both $D_s$  and $D_n$, $D_s\cap D_n$ has non-isomorphic
induced uniformizing systems from $(\tilde{D}_s,\Z_2,\pi_s)$ and
$(\tilde{D}_n,\Z_3,\pi_n)$, although they define the same germ at
each point in $D_s\cap D_n$. This also shows that although both
$D_s$ and $D_n$ are uniformized, their union $S^2$ can not be
uniformized, therefore is not a global quotient.} \vskip 0.1in

The notion of orbifold was first introduced by Satake in \cite{S},
where a different name, $V$-manifold, was used. Satake's
V-manifold corresponds to reduced orbifold in our case.  In [S],
an orbifold structure on a topological space $X$ is given by an
open cover $\U$ of $X$ satisfying the following conditions:
\begin{itemize}
\item {(2.1.1a)} Each element $U$ in $\U$ is uniformized, say by $(V,G,\pi)$.
\item {(2.1.1b)} If $U^\prime\subset U$, then there is a collection of
injections $(V^\prime,G^\prime,\pi^\prime)\rightarrow (V,G,\pi)$.
\item {(2.1.1c)} For any point $p\in U_1\cap U_2$, $U_1,U_2\in\U$,
there is a $U_3\in\U$ such that $p\in U_3\subset U_1\cap U_2$.
\end{itemize}
One can show that our definition is equivalent to Satake's.

Next we consider a class of continuous maps between two orbifolds
which carry an additional structure of differentiability with
respect to the orbifold structures. Let $U$ be uniformized by
$(V,G,\pi)$ and $U^\prime$ by $(V^\prime,G^\prime,\pi^\prime)$,
and $f: U\rightarrow U^\prime$ be a continuous map. A {\it $C^l$
lifting, $0\leq l\leq\infty$}, of $f$ is a $C^l$ map
$\tilde{f}:V\rightarrow V^\prime$ and a homomorphism $\lambda:
G\rightarrow G^\prime$ such that
$\pi^\prime\circ\tilde{f}=f\circ\pi$, and
$\lambda(g)\cdot\tilde{f}(x)=\tilde{f}(g\cdot x)$ for any $x\in
V$. Two liftings $\tilde{f}_i: (V_i,G_i,\pi_i)\rightarrow
(V^\prime_i,G^\prime_i,\pi^\prime_i)$, $i=1,2$, are {\it
isomorphic} if there exist isomorphisms
$(\phi,\tau):(V_1,G_1,\pi_1)\rightarrow (V_2,G_2,\pi_2)$ and
$(\phi^\prime,\tau^\prime):(V_1^\prime,G_1^\prime,\pi_1^\prime)\rightarrow
(V_2^\prime,G_2^\prime,\pi_2^\prime)$ such that $\phi^\prime\circ
\tilde{f}_1= \tilde{f}_2\circ \phi$. Let $p\in U$ for any
uniformized neighborhood $U_p$ of $p$ and uniformized neighborhood
$U_{f(p)}$ of $f(p)$ such that $f(U_p)\subset U_{f(p)}$, a lifting
$\tilde{f}$ of $f$ will induce a lifting $\tilde{f}_p$ for
$f|_{U_{p}}:U_p\rightarrow U_{f(p)}$ as follows: For any injection
$(\phi,\tau):(V_p,G_p,\pi_p)\rightarrow (V,G,\pi)$, consider the
map $\tilde{f}\circ\phi: V_p\rightarrow V^\prime$, and observe
that $\pi^\prime\circ\tilde{f}\circ\phi(V_p)\subset U_{f(p)}$
implies $\tilde{f}\circ\phi(V_p)\subset
(\pi^\prime)^{-1}(U_{f(p)})$. Therefore there is an injection
$(\phi^\prime,\tau^\prime):(V_{f(p)}, G_{f(p)},\pi_{f(p)})
\rightarrow (V^\prime, G^\prime, \pi^\prime)$ such that
$\tilde{f}\circ\phi(V_p)\subset \phi^\prime(V_{f(p)})$. We define
$\tilde{f}_p=(\phi^\prime)^{-1}\circ\tilde{f}\circ\phi$. In this
way we obtain a lifting $\tilde{f}_p: (V_p,G_p,\pi_p)\rightarrow
(V_{f(p)}, G_{f(p)},\pi_{f(p)})$ for $f|_{U_{p}}:U_p\rightarrow
U_{f(p)}$. We can verify that different choices give isomorphic
liftings. We define the {\it germ} of liftings as follows: two
liftings are {\it equivalent at $p$} if they induce isomorphic
liftings on a smaller neighborhood of $p$.

Now consider orbifolds $X$ and $X^\prime$ and a continuous map
$f:X\rightarrow X^\prime$. A {\it lifting} of $f$ consists of the
following data: for any point $p \in X$, there exist charts
$(V_p,G_p,\pi_p)$ at $p$ and $(V_{f(p)},G_{f(p)}, \pi_{f(p)})$ at
$f(p)$ and a lifting $\tilde{f}_p$ of $f_{\pi_p(V_p)}:
\pi_p(V_p)\rightarrow \pi_{f(p)}(V_{f(p)})$ such that for any
$q\in \pi_p(V_p)$, $\tilde{f}_p$ and $\tilde{f}_q$ induce the same
germ of liftings of $f$ at $q$. We can define the {\it germ} of
liftings in the sense that two liftings of $f$ $\{\tilde{f}_{p,i}:
(V_{p,i},G_{p,i},\pi_{p,i}) \rightarrow
(V_{f(p),i},G_{f(p),i},\pi_{f(p),i}): p\in X\}$, $i=1,2$, are {\it
equivalent} if for each $p\in X$, $\tilde{f}_{p,i}, i=1,2$, induce
the same germ of liftings of $f$ at $p$.

    \vskip 0.1in
\noindent{\bf Definition 2.1.6: }{\it A {\it $C^l$ map} ($0\leq
l\leq\infty$) between orbifolds $X$ and $X^\prime$ is a germ of
$C^l$ liftings of a continuous map between $X$ and $X^\prime$. We
denote by $\tilde{f}$ a $C^l$ map which is a germ of liftings of a
continuous map $f$.

A sequence of $C^l$ maps $\tilde{f}_n$ is said to converge to a
$C^l$ map $\tilde{f}_0$ in the $C^l$ topology if there exists a
sequence of liftings $\tilde{f}_{p,n}: (V_p,G_p,\pi_p)\rightarrow
(V_{f_n(p)},G_{f_n(p)},\pi_{f_n(p)})$ defining the germs
$\tilde{f}_n$ such that for each $p\in X$ there exists a chart
$(V_{f_0(p)},G_{f_0(p)},\pi_{f_0(p)})$ and an integer $n(p)>0$
with the following property: for each $n\geq n(p)$, there is an
injection $(\psi_{p,n},\tau_{p,n}):
(V_{f_n(p)},G_{f_n(p)},\pi_{f_n(p)})\rightarrow
(V_{f_0(p)},G_{f_0(p)},\pi_{f_0(p)})$ such that
$\psi_{p,n}\circ\tilde{f}_{p,n}$ converges in $C^l$ to
$\tilde{f}_{p,0}$ which defines the germ $\tilde{f}_0$.}\vskip
0.1in

\noindent{\bf Example 2.1.7a: }{\it The real line $\R$ as a smooth
manifold is trivially an orbifold. A $C^l$ map from an orbifold
$X$ to $\R$ is called a {\it $C^l$ function} on $X$. The set of
all $C^l$ functions on $X$ is denoted by $C^l(X)$. }\vskip 0.1in

\noindent{\bf Example 2.1.7b: }{\it Let $X=\R\times \C$, and be
given an orbifold structure by $(\R\times\C,\Z_4, \pi)$ where
$\Z_4$ acts only on the factor $\C$ by multiplication by
$\sqrt{-1}$. Define $C^1$ maps $\tilde{f}_1:\R\rightarrow
(\R\times\C,\Z_4,\pi)$ by $t\rightarrow (t,t^2)$ and
$\tilde{f}_2:\R\rightarrow (\R\times\C,\Z_4,\pi)$ by $t\rightarrow
(t,t^2)$ for $t\leq 0$ and $(t,\sqrt{-1}t^2)$ for $t\geq 0$. Then
$\tilde{f}_1,\tilde{f}_2$ induce the same continuous map
$f:\R\rightarrow X$, but they are {\it not} isomorphic as $C^1$
maps.}\vskip 0.1in

Next we describe the notion of orbifold vector bundle which
corresponds to the notion of smooth vector bundle. When there is
no confusion, we simply call it a vector bundle.  We begin with
local uniformizing systems for vector bundles. Given a uniformized
topological space $U$ and a topological space $E$ with a
surjective continuous map $pr:E\rightarrow U$, a {\it uniformizing
system of rank $k$ vector bundle} for $E$ over $U$ consists of the
following data:

\begin{itemize}
\item A uniformizing system $(V,G,\pi)$ of $U$.
\item A uniformizing system $(V\times\R^k,G,\tilde{\pi})$ for $E$.
The action of $G$ on $V\times\R^k$ is an extension of the action
of $G$ on $V$ given by $g(x,v)=(gx,\rho(x,g)v)$, where
$\rho:V\times G\rightarrow Aut(\R^k)$ is a smooth map satisfying:
$$ \rho(gx,h)\circ\rho(x,g)=\rho(x,h\circ g), \hspace{3mm} g,h\in
G, x\in V. $$
\item The natural projection map $\tilde{pr}:V\times\R^k\rightarrow V$
satisfies $\pi\circ\tilde{pr}=pr\circ\tilde{\pi}$.
\end{itemize}

We can similarly define {\it isomorphisms} between uniformizing
systems of  vector bundle for $E$ over $U$. The only additional
requirement is that the diffeomorphisms between $V\times\R^k$ are
linear on each fiber of $\tilde{pr}:V\times\R^k\rightarrow V$.
Moreover, for each connected open subset $U^\prime$ of $U$, we can
similarly prove that there is a unique isomorphism class of
induced uniformizing systems of vector bundle for
$E^\prime=pr^{-1}(U^\prime)$ over $U^\prime$. The {\it germ} of
uniformizing systems of vector bundle at a point $p\in U$ can also
be  similarly defined.

    \vskip 0.1in
\noindent{\bf Definition 2.1.8: }{\it Let $X$ be an orbifold and
$E$ be a topological space with a surjective continuous map
$pr:E\rightarrow X$. A {\it rank $k$ vector bundle  structure} on
$E$ over $X$ consists of the following data: For each point $p\in
X$, there is a uniformized neighborhood $U_p$ and a uniformizing
system of rank $k$  vector bundle for $pr^{-1}(U_p)$ over $U_p$
such that for any $q\in U_p$, the uniformizing systems of vector
bundle over $U_p$ and $U_q$ define the same germ at $q$. The {\it
germ} of rank $k$ vector bundle structures on $E$ over $X$ can be
similarly defined. The topological space $E$ with a given germ of
vector bundle structures becomes an orbifold and is called a {\it
vector bundle} over $X$. Each chart
$(V_p\times\R^k,G_p,\tilde{\pi}_p)$ is called a {\it local
trivialization} of $E$. At each point $p\in X$, the fiber
$E_p=pr^{-1}(p)$ is isomorphic to $\R^k/G_p$. It contains a linear
subspace $E^p$ of fixed points of $G_p$. Two vector bundles
$pr_1:E_1\rightarrow X$ and $pr_2:E_2\rightarrow X$ are {\it
isomorphic} if there is a $C^\infty$ map
$\tilde{\psi}:E_1\rightarrow E_2$ given by $\tilde{\psi}_p:
(V_{1,p}\times\R^k,G_{1,p},\tilde{\pi}_{1,p})\rightarrow
(V_{2,p}\times\R^k,G_{2,p},\tilde{\pi}_{2,p})$ which induces an
isomorphism between $(V_{1,p},G_{1,p},\pi_{1,p})$ and
$(V_{2,p},G_{2,p},\pi_{2,p})$, and is a linear isomorphism between
the fibers of $\tilde{pr}_{1,p}$ and $\tilde{pr}_{2,p}$. By
replacing $\R^k$ with $\C^k$, we have the definition of  complex
vector bundle}. \vskip 0.1in

\noindent{\bf Remark 2.1.10: }{\it There is a notion of vector
bundle over an orbifold with boundary. One can easily verify that
if $pr:E\rightarrow X$ is a  vector bundle over a orbifold with
boundary $X$, then the restriction to the boundary $\partial X$,
$E_{\partial X}=pr^{-1}(\partial X)$, is a vector bundle over
$\partial X$.}\vskip 0.1in

\noindent{\bf Remark 2.1.11: }{\it One can define orbifold fiber
bundle with fiber a general space in the same manner.}
    \vskip 0.1in

A $C^l$ map $\tilde{s}$ from $X$ to a  vector bundle
$pr:E\rightarrow X$ is called a {\it $C^l$ section} if locally
$\tilde{s}$ is given by $\tilde{s}_p: V_p\rightarrow
V_p\times\R^k$ where $\tilde{s}_p$ is $G_p$-equivariant and
$\tilde{pr}\circ\tilde{s}_p=Id$ on $V_p$. We observe that

\begin{itemize}
\item For each point $p$, $s(p)$ lies in $E^p$, the linear subspace of fixed
points of $G_p$.
\item The space of all $C^l$ sections of $E$, denoted by $C^l(E)$, has a
structure of vector space over $\R$ (or $\C$) as well as a
$C^l(X)$-module structure.
\item The $C^l$ sections $\tilde{s}$ are in $1:1$ correspondence with the
underlying continuous maps $s$.
\end{itemize}

\vskip 0.1in
    \noindent
    {\bf Remark 2.1.12: }{\it If $E\rightarrow X$ is an orbifold vector
    bundle over $X$ which does not introduce  an orbifold bundle over the associated
    reduced orbifold, then $E$ has no
    nonzero section. }
    \vskip 0.1in

     Orbifold vector bundles are more conveniently described by
transition maps (see [S]). More precisely, an orbifold vector
bundle over an orbifold $X$ can be constructed from the following
data: A compatible cover $\U$ of $X$ such that for any injection
$i:(V^\prime,G^\prime,\pi^\prime)\rightarrow (V,G,\pi)$, there is
a smooth map $g_i:V^\prime\rightarrow Aut(\R^k)$ giving an open
embedding $V^\prime\times\R^k\rightarrow V\times\R^k$ by
$(x,v)\rightarrow (i(x),g_i(x)v)$, and for any composition of
injections $j\circ i$, we have $$ g_{j\circ i}(x)=g_j(i(x))\circ
g_i(x), \forall x\in V. \leqno (2.1.2) $$ Two collections of maps
$g^{(1)}$ and $g^{(2)}$ define isomorphic bundles if there are
maps $\delta_V:V\rightarrow Aut(\R^k)$ such that for any injection
$i: (V^\prime,G^\prime,\pi^\prime) \rightarrow (V,G,\pi)$, we have
$$ g^{(2)}_i(x)=\delta_V(i(x))\circ g^{(1)}_i(x)\circ
(\delta_{V^\prime}(x))^{-1}, \forall x\in V^\prime. \leqno (2.1.3)
$$ Since $(2.1.2)$ behaves naturally under constructions of vector
spaces such as tensor product, exterior product, etc. we can
define these constructions for vector bundles.

    \vskip 0.1in
\noindent{\bf Example 2.1.13: }{\it For an orbifold $X$, the
tangent bundle $TX$ can be constructed because the differential of
any injection satisfies $(2.1.2)$. Likewise, we define the
cotangent bundle $T^\ast X$, and the bundles of exterior power or
tensor product. All the differential geometry such as de Rham
theory, connections, curvature and characteristic classes extends
to orbifold vector bundles.  Moreover, the de Rham cohomology of
an orbifold is isomorphic to the de Rham cohomology of its
associated reduced  orbifold, which is isomorphic to the singular
cohomology of the underlying topological space. Observe also that
if $\omega$ is a differential form on $X^\prime$ and
$\tilde{f}:X\rightarrow X^\prime$ is a $C^\infty$ map, then there
is a pull-back form $\tilde{f}^\ast\omega$ on $X$.} \vskip 0.1in

\subsection{Pull-back bundles and good maps}

Let $pr: E\rightarrow Y$ be a vector bundle over a topological
space $Y$. Then for any continuous map $f: X\rightarrow Y$ from  a
topological space $X$, the pull-back vector bundle $f^\ast E$ over
$X$ is well-defined. However, this is no longer the case for
orbifold vector  bundles. Let $pr:E\rightarrow X^\prime$ be an
orbifold vector bundle over $X^\prime$, and
$\tilde{f}:X\rightarrow X^\prime$ a $C^\infty$ map. By a {\it
pull-back bundle of $E$ over $X$ via $\tilde{f}$ } we mean an
orbifold vector bundle $\pi:E^\bullet\rightarrow X$ together with
a $C^\infty$ map $\bar{f}:E^\bullet\rightarrow E$ such that each
local lifting of $\bar{f}$ is an isomorphism restricted to each
fiber, and $\bar{f}$ covers the $C^\infty$ map $\tilde{f}$ between
the bases.

Let $\tilde{f}:X\rightarrow X^\prime$ be a $C^\infty$ map between
orbifolds $X$ and $X^\prime$ whose underlying continuous map is
denoted by $f$. Suppose there is a compatible cover $\U$ of $X$,
and a collection of open subsets $\U^\prime$ of $X^\prime$
satisfying $(2.1.1a-c)$ and the following condition: There is a
1:1 correspondence between elements of $\U$ and $\U^\prime$, say
$U\leftrightarrow U^\prime$, such that $f(U)\subset U^\prime$, and
an inclusion $U_2\subset U_1$ implies an inclusion
$U_2^\prime\subset U_1^\prime$. Moreover, there is a collection of
local $C^\infty$ liftings  $\{\tilde{f}_{UU^\prime}\}$ of $f$,
where $\tilde{f}_{UU^\prime}:(V,G,\pi)\rightarrow
(V^\prime,G^\prime,\pi^\prime)$ satisfies the following condition:
each  injection $i:(V_2,G_2,\pi_2)\rightarrow (V_1,G_1,\pi_1)$ is
assigned  an injection
$\lambda(i):(V_2^\prime,G_2^\prime,\pi^\prime_2)\rightarrow
(V^\prime_1,G_1^\prime,\pi_1^\prime)$ such that
$\tilde{f}_{U_1U^\prime_1}\circ
i=\lambda(i)\circ\tilde{f}_{U_2U^\prime_2}$, and for any
composition of injections $j\circ i$, the following compatibility
condition holds: $$ \lambda(j\circ i)=\lambda(j)\circ \lambda(i).
\leqno (2.2.1) $$ Observe that when the injection
$i:(V,G,\pi)\rightarrow (V,G,\pi)$ is an automorphism of
$(V,G,\pi)$, the assignment of $\lambda(i)
:(V^\prime,G^\prime,\pi^\prime)\rightarrow
(V^\prime,G^\prime,\pi^\prime)$ to $i$ satisfying $(2.2.1)$ is
equivalent to a homomorphism $\lambda_{UU^\prime}: G\rightarrow
G^\prime$. We call $\lambda_{UU^\prime}: G\rightarrow G^\prime$
the {\it group homomorphism} of
$\{\tilde{f}_{UU^\prime},\lambda\}$ on $U$. Such a collection of
maps clearly defines a $C^\infty$ lifting of the continuous map
$f$. If it is in the same germ as $\tilde{f}$, we call
$\{\tilde{f}_{UU^\prime},\lambda\}$ a {\it compatible system} of
$\tilde{f}$.

    \vskip 0.1in
\noindent{\bf Definition 2.2.1: }{\it A $C^\infty$ map is called
{\it good} if it admits a compatible system.}\vskip 0.1in

\noindent{\bf Example 2.2.2: }{\it  There can be essentially
different compatible systems of the same $C^\infty$ map, as shown
in the following example: Let $X=\C\times\C/G$ where $G=\Z_2
\oplus\Z_2$ acting on $\C\times\C$ in the standard way. For the
$C^l$ map $\tilde{f}:(\C,\Z_2)\rightarrow (\C\times\C,G)$ defined
by the inclusion of $\C\times\{0\}$, there are two compatible
systems $(\tilde{f},\lambda_i): (\C,\Z_2)\rightarrow
(\C\times\C,G)$, $i=1,2$, for $\lambda_1(1)=(1,0)$ and
$\lambda_2(1)=(1,1)$, which are apparently different.}\vskip 0.1in

\noindent{\bf Lemma 2.2.3: }{\it Let $pr:E\rightarrow X^\prime$ be
an orbifold vector bundle over $X^\prime$. For any $C^\infty$
compatible system $\xi=\{\tilde{f}_{UU^\prime},\lambda\}$ of a
good $C^\infty$ map $\tilde{f}:X\rightarrow X^\prime$, there is a
canonically constructed pull-back bundle of $E$ via $\tilde{f}$: a
 bundle $pr:E^\bullet_\xi\rightarrow X$ together with a $C^\infty$ map
$\bar{f}_\xi:E^\bullet_\xi \rightarrow E$ covering $\tilde{f}$.
Let $c$ be a universal characteristic class defined by the
Chern-Weil construction; then
$\tilde{f}^\ast(c(E))=c(E^\bullet_\xi)$.}\vskip 0.1in

\noindent{\bf Proof:} Without loss of generality, we assume that
$E_{U^\prime}$ over $U^\prime$ is uniformized by $(V^\prime\times
\R^k,G^\prime,\tilde{\pi}^\prime)$. Then we have a collection of
pull-back bundles $\tilde{f}_{UU^\prime}^\ast (V^\prime \times
\R^k)$ over $V$ which have the form of $V\times\R^k$. Let
$\{g^\prime\}$ be a collection of transition maps of $E$ with
respect to $\U^\prime$; we define a set of transition maps $\{g\}$
on $X$ with respect to the cover $\U$ by pull-backs, i.e., we set
$g_i=g^\prime_{\lambda(i)}\circ\tilde{f}_{U_2U^\prime_2}$ for any
injection $i:(V_2,G_2,\pi_2)\rightarrow (V_1,G_1,\pi_1)$, where
$\lambda(i):(V_2^\prime,G_2^\prime,\pi_2^\prime)\rightarrow
(V_1^\prime,G_1^\prime,\pi_1^\prime)$ is the injection assigned to
$i$. Then the compatibility condition $(2.1.1)$ implies that the
set of maps $\{g\}$ satisfies  equation $(2.1.2)$, which defines a
bundle over $X$. We denote it by $pr:E^\bullet\rightarrow X$. The
existence of a $C^\infty$ map $\bar{f}:E^\bullet\rightarrow E$ is
obvious from the construction. On the other hand, for any
connection $\nabla$ on $E$, there is a pull-back connection
$\tilde{f}^\ast(\nabla)$ on $E^\bullet$, so that the equation
$\tilde{f}^\ast(c(E))=c(E^\bullet)$ holds for any universal
characteristic class $c$ defined by the Chern-Weil construction.
\hfill $\Box$

    \vskip 0.1in
\noindent{\bf Definition 2.2.4: }{\it  Two compatible systems
$\xi_i$ for $i=1,2,$ of a $C^\infty$ map $\tilde{f}:X\rightarrow
X^\prime$ are {\it isomorphic} if for any orbifold vector bundle
$E$ over $X^\prime$ there is an isomorphism $\psi$ between the
corresponding pull-back bundles $E^\ast_i$ with
$\bar{f}_i:E_i^\ast\rightarrow E$, $i=1,2$, such that
$\bar{f}_1=\psi\circ\bar{f}_2$.}\vskip 0.1in

    \vskip 0.1in
\noindent{\bf Definition 2.2.5: }{\it  It is easily seen that for
each $p\in X$, a compatible system determines a group homomorphism
$G_p\rightarrow G_{f(p)}$, and an isomorphism class of compatible
systems determines a conjugacy class of homomorphisms. We call
such a conjugacy class of group homomorphisms the {\it group
homomorphism} of the said isomorphism class of compatible systems
at $p$.} \vskip 0.1in

    Following is a very important class of examples of good maps.

\noindent{\bf Definition 2.2.6: }{\it A $C^\infty$ map
$\tilde{f}:X\rightarrow X^\prime$ between orbifolds is called {\it
regular} if the underlying continuous map $f$ has the following
property: $f^{-1}(X_{reg}^\prime)$ is an open dense and connected
subset of $X$.}\vskip 0.1in

    The following important lemma is proved by Chen-Ruan \cite{CR3}.

\noindent{\bf Lemma 2.2.7: }{\it If $\tilde{f}$ is regular, then
$\tilde{f}$ is the unique germ of $C^\infty$ liftings of $f$.
Moreover, $\tilde{f}$ is good with a unique isomorphism class of
compatible systems.} \vskip 0.1in

\noindent{\bf Remark 2.2.8: }{\it  Here are some examples of
regular $C^\infty$ maps. Let $E$ be a bundle over a reduced
orbifold $X$. Then any $C^\infty$ section $\tilde{s}$ of $E$, as a
$C^\infty$ map $X\rightarrow E$, is regular, since $s(p)$ is in
$\Sigma E$ only if $p$ is in $\Sigma X$, which is of codimension
at least two. Another example: Let $E^\bullet$ be a pull-back
bundle over $X$ of the tangent bundle of $TX^\prime$ with the
$C^\infty$ map $\bar{f}:E^\bullet\rightarrow TX^\prime$; then
$\bar{f}$ is a regular map, since if $\bar{f}(p,v)$ is in $\Sigma
TX^\prime$, then $f(p)$ is in $\Sigma X^\prime$ and $v$ is in a
subset of $(TX^\prime)_{f(p)}$ of codimension at least two.}\vskip
0.1in

\noindent{\bf Remark 2.2.9: }{\it  Here is another class of
regular $C^\infty$ maps. A $C^\infty$ map $\tilde{f}: X\rightarrow
X^\prime$ is called a {\it $C^\infty$ embedding} if each local
lifting $\tilde{f}_p:(V_p,G_p,\pi_p)\rightarrow
(V_{f(p)},G_{f(p)},\pi_{f(p)})$ is a $\lambda_p$-equivariant
embedding for some isomorphism $\lambda_p: G_p\rightarrow
G_{f(p)}$. It is easily seen that $f^{-1}(\Sigma X^\prime) =\Sigma
X$, so that $\tilde{f}$ is regular. For a $C^\infty$ embedding
$\tilde{f}:X\rightarrow X^\prime$, the normal bundle of $f(X)$ in
$X^\prime$ as a bundle is well-defined, and is isomorphic to the
quotient bundle $(TX^\prime)^\bullet/TX$.}\vskip 0.1in

\noindent{\bf Example 2.2.10: }{\it  Let $\P^1=\{[z_0,z_1]\}$ be
the 1-dimensional complex projective space. We define a $\Z_2$
action on it by $x\cdot [z_0,z_1]=[xz_0,z_1]$. Let $X=\P^1/\Z_2$
be the orbifold as quotient space. Similarly, let $\P^2
=\{[z_0,z_1,z_2]\}$ be the 2-dimensional complex projective space,
with a $\Z_2\oplus \Z_2$ action on it, given by $(x,y)\cdot
[z_0,z_1,z_2] =[xz_0,yz_1,z_2]$. Let $X^\prime=\P^2/(\Z_2\oplus
\Z_2)$ be the orbifold as quotient space. We consider two
sequences of $C^\infty$ maps (actually they are holomorphic)
$\tilde{f}_n,\tilde{g}_n:X\rightarrow X^\prime$ defined by
$\tilde{f}_n([z_0,z_1])=[z_0,n^{-1}z_1,z_1]$ and
$\tilde{g}_n([z_0,z_1])=[z_0,n^{-1}z_0,z_1]$. It is easily seen
that both sequences consist of regular maps, so they are good
maps. As $n\rightarrow \infty$, both sequences converge. Let
$\tilde{f}= \lim\tilde{f}_n$ and $\tilde{g}=\lim\tilde{g}_n$. Then
both $\tilde{f}$ and $\tilde{g}$ are good maps (as we shall see),
and $\tilde{f}=\tilde{g}$ as $C^\infty$ maps, but with different
isomorphism classes of compatible systems. In fact, the group
homomorphism of $\tilde{f}$ which is from $\Z_2$ to $\Z_2\oplus
\Z_2$ is  given by $x\rightarrow (x,1)$ at $[0,1]$, and
$x\rightarrow (x,x)$ at $[1,0]$. For $\tilde{g}$, it is given by
$x\rightarrow (x,x)$ at $[0,1]$ and $x\rightarrow (1,x)$ at
$[1,0]$.}\vskip 0.1in \hfill $\Box$

The operation of composition is well-defined for good maps, as
shown in the following lemma.

    \vskip 0.1in
\noindent{\bf Lemma 2.2.11: }{\it  Let $\tilde{f}$, $\tilde{g}$ be
two good $C^\infty$ maps; then the composition
$\tilde{g}\circ\tilde{f}$ is also a good $C^\infty$ map, and any
isomorphism class of compatible systems of $\tilde{f}$ and
$\tilde{g}$ determines a unique isomorphism class of compatible
systems for the composition $\tilde{g}\circ\tilde{f}$.}

\vspace{2mm}

Now consider a good $C^\infty$ map $\tilde{f}:X\rightarrow
X^\prime$ with an isomorphism class of compatible systems $\xi$.
Then we have an isomorphism class of pull-back bundles
$(TX^\prime)^\bullet_\xi$ over $X$ and the good $C^\infty$ map
$\bar{f}_\xi:(TX^\prime)^\bullet_\xi\rightarrow TX^\prime$. For
any $C^\infty$ section $\tilde{s}$ of $(TX^\prime)^\bullet_\xi$,
we take the composition
$\tilde{f}_{\xi,s}=Exp\circ\bar{f}_\xi\circ\tilde{s}$ from $X$
into $X^\prime$. Then $\tilde{f}_{\xi,s}$ is a good $C^\infty$ map
with an isomorphism class of compatible systems determined by
$\xi$. A natural question is: given a good map $\tilde{g}$ nearby
$\tilde{f}$ with an isomorphism class of compatible systems, is
there an isomorphism class of compatible systems $\xi$ of
$\tilde{f}$, and a $C^\infty$ section $\tilde{s}$ of the pull-back
bundle $(TX^\prime)^\bullet_{\xi}$ such that $\tilde{g}$ is
realized as $\tilde{f}_{\xi,s}$ ?  If there is, will $\xi$ and
$\tilde{s}$ be unique? These questions seem to be non-trivial in
general, but as we shall see, it can be dealt with in certain
special cases, e.g., when $\tilde{f}$, $\tilde{g}$ are
pseudo-holomorphic maps from a complex orbicurve into an almost
complex orbifold. We refer reader to \cite{CR3} for details.

\section{Orbifold Cohomology}
   As I mentioned in the introduction, the ordinary cohomology is
   a wrong theory for orbifolds. The correct one has to incorporate
   the twisted sectors. Furthermore, the internal
   freedom of orbifold string theory allows a twisting as well.
   Such a theory (orbifold cohomology) without twisting was constructed
   by Chen-Ruan \cite{CR1}. The twisted orbifold cohomology was constructed
   by Ruan \cite{R}. In the case of global quotients, orbifold cohomology
   groups
   were known to physicists  \cite{Z},
   \cite{VW}. However, even in this case, the orbifold cup product is new. This section is a combination of \cite{CR1},
   \cite{R}.

    \subsection{Twisted sector and inner local system}

Let $X$ be an orbifold. For any point $p\in X$, let
$(V_p,G_p,\pi_p)$ be a local chart at $p$. Let
$\widetilde{\Sigma_k X}$ denote the set of pairs $(p,(\g))$, where
$(\g)$ stands for the conjugacy class of $\g=(g_1,\cdots,g_k)$ by
an element of $G_p$.  We call the $\widetilde{X}_k$ {\em
multi-sectors}.

\vspace{3mm}

\noindent{\bf Lemma 3.1.1:} {\it The multi-sector
$\widetilde{\Sigma_k X}$ is naturally an orbifold,  and is a
finite union of closed orbifolds when $X$ is closed, with the
orbifold structure given by $$
\{\pi_{p,\g}:(V^\g_p,C(\g))\rightarrow V^\g_p/C(\g); (p,(\g))\in
\widetilde{\Sigma_k X}\} \leqno(3.1.1) $$ where
$V^\g_p=V_p^{g_1}\cap V_p^{g_2}\cap\cdots\cap V_p^{g_k}$,
$C(\g)=C(g_1)\cap C(g_2)\cap\cdots\cap C(g_k)$. Here
$\g=(g_1,\cdots,g_k)$, $V_p^g$ stands for the fixed-point set of
$g\in G_p$ in $V_p$, and $C(g)$ for the centralizer of $g$ in
$G_p$,  for some local chart $(V_p,G_p,\pi_p)$ at $p$.}

\vspace{3mm}

\noindent{\bf Proof:} First we identify a point $(q,(\h))$ in
$\widetilde{\Sigma_k X}_k$ as a point in $\bigsqcup_{\{(p,(\g))\in
\widetilde{\Sigma_k X}\}}V_p^{\g}/C(\g)$ if $q\in U_p$. Pick a
representative $y\in V_p$ such that $\pi_p(y)=q$. Then this gives
rise to a monomorphism $\lambda_y:G_q\rightarrow G_p$. Pick a
representative $\h=(h_1,\cdots,h_k)\in G_q\times \cdots\times G_q$
for $(\h)$, and let $\g=\lambda_y(\h)$. Then $y\in V^{\g}_p$. So
we have a map $\theta:(q,\h)\rightarrow (y,\g)$. If we change $\h$
by $\h^\prime=a^{-1}\h a$ for some $a\in G_q$, then $\g$ is
changed to $\lambda_y(a^{-1}\h a)=\lambda(a)^{-1}\g\lambda(a)$. So
we have $\theta:(q,a^{-1}\h a)\rightarrow
(y,\lambda(a)^{-1}\g\lambda(a))$ where $y$ is regarded as an
element in $V^{\lambda(a)^{-1}\g\lambda(a)}_p$. (Note that
$\lambda_y$ is determined up to conjugacy by an element in $G_q$.)
If we take a different representative $y^\prime\in V_p$ such that
$\pi_p(y^\prime)=q$, and assume $y^\prime=b\cdot y$ for some $b\in
G_p$. then we have a different identification
$\lambda_{y^\prime}:G_q\rightarrow G_p$ of $G_q$ as a subgroup of
$G_p$ where $\lambda_{y^\prime}=b\cdot\lambda_y\cdot b^{-1}$. In
this case, we have $\theta:(q,\h)\rightarrow (y^\prime,b\g
b^{-1})$ where $y^\prime\in V_p^{b\g b^{-1}}$. If $\g=b\g b^{-1}$,
then $b\in C(\g)$. Therefore we have shown that $\theta$ induces a
map sending $(q,(\h))$ to a point in
$\bigsqcup_{\{(p,(\g))\in\widetilde{\Sigma_k X}\}}V_p^{\g}/C(\g)$,
which can be similarly shown to be one to one and onto. Hence we
have shown that $\widetilde{\Sigma_k X}$ is covered by
$\bigsqcup_{\{(p,(\g))\in\widetilde{\Sigma_k X}\}}V_p^{\g}/C(\g)$.

We define a topology on $\widetilde{\Sigma_k X}$ so that each
$V_p^{\g}/C(\g)$ is an open subset for any $(p,\g)$. We also
uniformize $V_p^{\g}/C(\g)$ by $(V_p^{\g},C(\g))$.  It remains to
show that these charts fit together to form an orbifold structure
on $\widetilde{\Sigma_k X}$. Let $x\in V_p^{\g}/C(\g)$ and take a
representative $\tilde{x}$ in $V_p^\g$. Let $H_x$ be the isotropy
subgroup of $\tilde{x}$ in $C(\g)$. Then $(V_p^\g,C(\g))$ induces
a germ of uniformizing system at $x$ as $(B_x,H_x)$ where $B_x$ is
a small ball in $V_p^\g$ centered at $\tilde{x}$. Let
$\pi_p(\tilde{x})=q$. We need to write $(B_x,H_x)$ as
$(V_q^\h,C(\h))$ for some $\h\in G_q\times\cdots\times G_q$. We
let $\lambda_x: G_q\rightarrow G_p$ be the induced monomorphism
resulting from choosing $\tilde{x}$ as the representative of $q$
in $V_p$. We define $\h=(\lambda_x)^{-1}(\g)$ (each $g_i$ is in
$\lambda_x(G_q)$ since $\tilde{x}\in V^\g_p$ and
$\pi_p(\tilde{x})=q$.) Then we can identify $B_x$ as $V_q^\h$. We
also see that $H_x=\lambda_x(C(\h))$. Hence we have proved that
$\widetilde{\Sigma_k X}$ is naturally an orbifold with an orbifold
structure described above ($\widetilde{\Sigma_k X}$ is Hausdorff
and second countable with the given topology for similar reasons).
The rest of the lemma is obvious. \hfill $\Box$

\vspace{3mm}
    {\bf Remark 3.1.2:  }{\it $\widetilde{\Sigma X}$ was introduced by
    Kawasaki \cite{Ka} in relation to the index theorem. A connected
    component of $\widetilde{\Sigma X}$ is called a  sector
    and will contribute to the orbifold cohomology group.
    $\widetilde{\Sigma_2 X}$ will be used to construct the Poincar\'e
    pairing. $\widetilde{\Sigma_3 X}$ will be used to define cup product and
    $\widetilde{\Sigma_4 X}$ will be used to prove associativity
    of the
    orbifold product. $\widetilde{\Sigma_k X}$ corresponds to the higher product for
    $k\geq 4$.}
    \vskip 0.1in

    Next, we consider some natural maps between multi-sectors.
    There are evaluation maps
    $$e_{i_1, \cdots, i_l}: \widetilde{\Sigma_k X}\rightarrow \widetilde{\Sigma_l X}\leqno(3.1.2)$$
    defined by $e_{i_1, \cdots, i_l}(x, (\g))\rightarrow (x,
    (g_{i_1}, \cdots, g_{i_l}))$.
    There is an involution
    $$I: \widetilde{\Sigma_k X}\rightarrow \widetilde{\Sigma_k X}\leqno(3.1.3)$$
    defined by
    $I(x, (\g))=(x, (\g^{-1})$, where $\g^{-1}=(g^{-1}_1, \cdots,
    g^{-1}_k).$
    \vskip 0.1in
    \noindent
    {\bf Lemma 3.1.3 : }{\it $e_{i_1,\cdots, i_l}$ and $I$ are good maps.}
    \vskip 0.1in
    {\bf Proof: } Suppose that $\{(V_p, G_p, \pi_p)\}$ is the
    orbifold structure of $X$. By Lemma 3.1.3, it induces the orbifold
    structures $\{V^{\g}_p, C(\g))\}, \{ (V^{(g_{i_1}, \cdots, g_{i_l})}_p,
    C(g_{i_1}, \cdots, g_{i_l})\}$
    of $\widetilde{\Sigma_k X}, \widetilde{\Sigma_l X}$
    simultaneously. Together with obvious embedding, they give a
    compatible system for $e_{i_1, \cdots, i_l}$. The proof for
    $I$ is similar. We leave it to the reader. $\Box$

    Next, we would like
   to describe the connected components of $\widetilde{\Sigma_k X}$. Recall
   that every point $p$ has a local chart $(V_p,G_p,\pi_p)$ which gives a
   local uniformized neighborhood $U_p=\pi_p(V_p)$.
   If $q\in U_p$, up to conjugation, there is an
   injective homomorphism $G_q\rightarrow G_p$. For $\g\in G_q$,
   the conjugacy class $(\g)_{G_p}$ is well-defined. We define
   an equivalence relation $(\g)_{G_q}\cong (\g)_{G_p}$. Let $T_k$ be
   the set of equivalence classes.  With abuse of the notation, we
   often use $(\g)$ to denote the equivalence class which $(\g)_{G_q}$
   belongs to. Let $T^o_k\subset T_k$ be the equivalence class of
   $(\g)$ such that $g_1\cdots g_k=1$.

   It is clear that $\widetilde{\Sigma_k X}$ is decomposed as a disjoint
   union of connected components
   $$ \widetilde{\Sigma_k X}=\bigsqcup_{(\g)\in T_k} X_{(\g)},\leqno(3.1.4)$$
   where
   $$X_{(\g)}=\{(p,(\g')_{G_p})|\g'\in G_p, (\g')_{G_p}\in (g)\}.\leqno(3.1.5)$$
   \vskip 0.1in
   \noindent
   {\bf Definition 3.1.4: }{\it $X_{(g)}$ for $g\neq 1$ is called
   a twisted sector.
   Furthermore, we call $X_{(1)}=X$ the nontwisted sector.}
   \vskip 0.1in

\noindent{\bf Example 3.1.5:} Consider the case that the orbifold
$X=Y/G$ is a global quotient. We will show that $\widetilde{\Sigma
X}$ is given by $\bigsqcup_{\{(g),g\in G\}} Y^g/C(g)$ where $Y^g$
is the fixed-point set of element $g\in G$. Let
$\pi:\widetilde{\Sigma X}\rightarrow X$ be the surjective map
defined by $(p,(g))\rightarrow p$. Then for any $p\in X$, the
preimage $\pi^{-1}(p)$ in $\widetilde{\Sigma X}$ has a
neighborhood described by $W_p=\bigsqcup_{\{(g),g\in
G_p\}}V_p^g/C(g)$, which is uniformized by
$\widetilde{W}_p=\bigsqcup_{\{(g),g\in G_p\}} V_p^g$. For each
$p\in X$, pick a $y\in Y$ that represents $p$, and an injection
$(\phi_p,\lambda_p): (V_p,G_p)\rightarrow (Y,G)$ whose image is
centered at $y$. This induces an open embedding
$\tilde{f}_p:\widetilde{W}_p\rightarrow
\bigsqcup_{\{(\lambda_p(g)),\lambda_p(g)\in G\}}Y^{\lambda_p(g)}
\subset \bigsqcup_{\{(g),g\in G\}}Y^g$, which induces a
homeomorphism $f_p$ from $W_p$ into $\bigsqcup_{\{(g),g\in
G\}}Y^g/C(g)$ that is independent of the choice of $y$ and
$(\phi_p,\lambda_p)$. These maps $\{f_p;p\in X\}$ fit together to
define a map $f:\widetilde{\Sigma X}\rightarrow \bigsqcup_{\{(g),
g\in G\}}Y^g/C(g)$ which we can verify to be a homeomorphism.

\hfill $\Box$

\vspace{3mm}

     Now, we introduce the notion of inner local system for orbifold.
    \vskip 0.1in
    \noindent
    {\bf Definition 3.1.6: }{\it Suppose that $X$ is an orbifold (almost complex or
    not). An inner local system $\L=\{L_{(g)}\}_{g\in T_1}$ is an assignment of  a flat
    complex line orbifold-bundle
    $$L_{(g)}\rightarrow X_{(g)}$$
    to each sector
    $X_{(g)}$ satisfying the compatibility condition
    \begin{description}
    \item[(1)] $L_{(1)}=1$ is trivial.
    \item[(2)] $I^*L_{(g^{-1})}=L^{-1}_{(g)}.$
    \item[(3)] Over each $X_{(\g)}$ with $(\g)\in T^o_3$,
    $\otimes_i e^*_i L_{(g_i)}=1 $.
    \end{description}
    If $X$ is a complex orbifold, we assume that $L_{(g)}$ is holomorphic.}
    \vskip 0.1in
    \noindent
    {\bf Lemma 3.1.7 :}{\it Suppose that $\L$ is an inner local system.
    For any $X_{(\g)}$ with $(\g)\in T^o_k$ for $k\geq 4$,
    $$\otimes_i e^*_i L_{(g_i)}=1.\leqno(3.1.6)$$
    }
    \vskip 0.1in
    {\bf Proof: } Suppose that $\g=(g_1, \cdots, g_k)$. We can
    define a sequence of triple elements
    $$\h_1=(g_1, g_2, (g_1g_2)^{-1}), \h_2=(g_1g_2,g_3 (g_1g_2g_3)^{-1}),
    \cdots, \h_{k-2}=(g_1\cdots g_{k-2}, g_{k-1}, g_k).$$
    By the construction, $(\h_i)\in T^o_3$. Moreover, the
    evaluation map $\prod_i e_i$ factors through the evaluation map
    to $\prod_i X_{(\h_i)}$. Then the lemma follows from (3).
    $\Box$

    An important way to produce inner local systems is by discrete
    torsion.

    First, we recall the definition of orbifold fundamental group.
                \vskip 0.1in
        \noindent
        {\bf Definition 3.1.8: }{\it A smooth map $f:Y\rightarrow X$
        is an orbifold cover iff (1) each $p\in Y$ has a neighborhood
      $U_p/G_p$ such that the restriction of $f$ to $U_p/G_p$ is    isomorphic
      to a map $U_p/G_p\rightarrow U_p/\Gamma$ such that $G_p\subset
    \Gamma$ is a subgroup. (2) Each $q\in X$ has a neighborhood
    $U_q/G_q$ for which each component of $f^{-1}(U_q/G_q)$ is
    isomorphic to $U_q/\Gamma'$ such that $\Gamma'\subset G_q$ is a
    subgroup.           An
        orbifold universal cover $f:Y\rightarrow X$ of $X$ has the
        property: (i) $Y$ is connected; (ii)if $f': Y'\rightarrow X$ is an orbifold
        cover, then there exists an orbifold cover $h:
        Y\rightarrow Y'$ such that $f=f'\circ h$. If $Y$ exists, we call $Y$ the
        orbifold universal cover of $X$ and the group of deck translations the orbifold
        fundamental group $\pi^{orb}_1(X)$ of $X$.}
        \vskip 0.1in
        By Thurston \cite{T}, an orbifold universal cover exists.
         It is clear from the definition that the orbifold universal
        cover is unique. Suppose that $f: Y\rightarrow X$ is an orbifold universal
        cover. Then
        $$f: Y-f^{-1}(\Sigma X)\rightarrow X-\Sigma X\leqno(3.1.7)$$
        is an honest cover with $G=\pi^{orb}_1(X)$ as covering
        group, where $\Sigma$ is the singular locus of $X$. Therefore $X=Y/G$ and there is a surjective
        homomorphism
        $$p_f: \pi_1(X-\Sigma X)\rightarrow G.\leqno(3.1.8)$$
        In general, (3.1.7) is not a universal covering. Hence,
        $p_f$ is not an isomorphism.

        \vskip 0.1in
        \noindent
        {\bf Remark 3.1.9 : }{\it Suppose that $X=Z/G$ for an orbifold $Z$
     and $Y$ is the orbifold universal cover of $Z$. By the
        definition, $Y$ is an orbifold universal cover of $X$. It
        is clear that there is a short exact sequence
        $$1\rightarrow \pi_1(Z)\rightarrow \pi^{orb}(X)\rightarrow
        G\rightarrow 1.\leqno(3.1.9)$$}

        \vskip 0.1in
        \noindent
        {\bf Example 3.1.10: } Consider the Kummer surface $T^4/\tau$
        where $\tau$ is the involution
        $$\tau(e^{it_1}, e^{it_2}, e^{it_3}, e^{it_4})=
        (e^{-it_1}, e^{-it_2}, e^{-it_3}, e^{-it_4}).\leqno(3.1.10)$$
        The universal cover is $\R^4$. The group $G$ of deck
        translations is generated by translations $\lambda_i$ by an integral point
        and the involution
    $$\tau: (t_1, t_2, t_3, t_4)\rightarrow (-t_1,
        -t_2, -t_3, -t_4).$$
     It is easy to check that
        $$G=\{\lambda_i (i=1,2,3,4), \tau| \tau^2=1,
        \tau\lambda_i=\lambda^{-1}_i\tau,\}\leqno(3.1.11)$$
        where $\lambda_i$ represents translation and $\tau$ represents
          involution.
        \vskip 0.1in
        \noindent
        {\bf Example 3.1.11: } Let $T^6=\R^6/\Gamma$ where $\Gamma$ is the lattice
        of integral points. Consider $\Z^2_2$ acting on $T^6$ lifted
        to an action on $\R^6$ as
        $$\sigma_1(t_1, t_2, t_3, t_4, t_5, t_6)=(-t_1, -t_2, -t_3, -t_4, t_5,
        t_6)$$
        $$\sigma_2(t_1, t_2, t_3, t_4, t_5, t_6)=(-t_1, -t_2, t_3, t_4, -t_5,
        -t_6)$$
        $$\sigma_3(t_1, t_2, t_3, t_4, t_5, t_6)=(t_1, t_2, -t_3, -t_4, -t_5,
        -t_6).$$
        This  example was considered by Vafa-Witten \cite{VW}.
        The orbifold fundamental group
        $$\pi^{orb}_1(T^6/\Z^2_2)=\{\tau_i (1\leq i\leq 6), \sigma_j
        (1\leq j\leq 3)|$$
        $$\sigma^2_i=1, \sigma_1\tau_i=\tau^{-1}_i\sigma_1 (i\neq
        5,6), \sigma_2\tau_i=\tau^{-1}_i\sigma_2 (i\neq 3,4),
        \sigma_3\tau_i=\tau^{-1}_i\sigma_3 (i\neq
        1,2)\}.\leqno(3.1.12)$$
        \vskip 0.1in
        The following example was taken from \cite{SC}
        \vskip 0.1in
        \noindent
        {\bf Example 3.1.12: }Consider the orbifold Riemann surface
        $\Sigma_g$ of genus $g$ and $n$ orbifold points $\z=(x_1, \cdots, x_n)$
        with orders $k_1, \cdots, k_n$. Then,
        $$\pi_1^{orb}(\Sigma_g)=\{\lambda_i (i\leq 2g), \sigma_i (i\leq n)|
        \sigma_1\cdots \sigma_n\prod_i [\lambda_{2i-1}, \lambda_{2i}]=1,
        \sigma^{k_i}_i=1\},\leqno(3.7)$$
        where $\lambda_i$ are the generators of $\pi_1(\Sigma_g)$ and
        $\sigma_i$ are the generators of $\Sigma_g-\z$ represented by a loop
        around each orbifold point.

         Note that $\pi_1^{orb}(\Sigma_g)$ is just  $\pi_1(\Sigma_g-\z)$ modulo
         the relation $\sigma^{k_i}_i=1$. This suggests that one can first take the cover
        of $\Sigma_g-\z$ induced by $\pi^{orb}_1(\Sigma)$. The relation $\sigma^{k_i}_i=1$
        implies that the preimage of the punctured disc around $x_i$ is a punctured disc.
        Then we can fill in the center point to obtain the orbifold
        universal cover.
        \vskip 0.1in
        \noindent
        {\bf Definition 3.1.13: }{ We call an element $\alpha\in H^2(\pi^{orb}_1(X), U(1))$
        a discrete torsion of $X$. }
        \vskip 0.1in
        If $X=Z/G$ for a finite group $G$, by Remark 3.2, there is a surjective homomorphism
        $$\pi: \pi^{orb}_1(X)\rightarrow G.$$
        $\pi$ induces a homomorphism
        $$\pi^*: H^2(G, U(1))\rightarrow H^2(\pi^{orb}_1(X),
        U(1)).\leqno(3.1.14)$$
        Hence, an element of $H^2(G, U(1))$ induces a discrete
        torsion of $X$.

        There are many ways to define $H^2(G, U(1))$. The definition $H^2(G, U(1))=H^2(BG, U(1))$ is a very useful definition
    for computation since we can use algebro-topological machinery. However,
          we can also take the original definition in terms of cocycles.
    A 2-cocycle is a map $\alpha: G\times
         G\rightarrow U(1)$ satisfying
    $$\alpha_{g,1}=\alpha_{1,g}=1, \alpha_{g,hk}\alpha_{h,k}
         =\alpha_{g,h}\alpha_{gh,k}, \leqno(3.1.15)$$
         for any $g,h,k\in G$. We denote the set of two-cocycles by
         $Z^2(G, U(1))$.
    For any map $\rho: G\rightarrow U(1)$ with $\rho_1=1$, its coboundary
    is defined by the formula
    $$(\delta \rho)_{g,h}=\rho_g\rho_h\rho_{gh}^{-1}.\leqno(3.1.16)$$
    Let $B^2(G, U(1))$ be the set
    of coboundaries. Then, $H^2(G,U(1))=Z^2(G,U(1))/B^2(G,U(1))$.
    $H^2(G, U(1))$ naturally appears in many important places in mathematics.
    For example, it classifies the group extensions of $G$ by $U(1)$. If
    we have a unitary projective representation of $G$, it naturally induces
    a class in $H^2(G, U(1))$. In many instances, this class completely classifies
    the projective unitary representation. In fact, it is in this context that
    discrete torsion arises in orbifold string theory.

    \vskip 0.1in
         \noindent
         {\bf Definition 3.1.13: }{\it
         For each 2-cocycle $\alpha$, we define its phase
         $$\gamma(\alpha)_{g,h}=\alpha_{g,h}\alpha^{-1}_{h,g}.\leqno(3.1.17)$$
         }
         \vskip 0.1in
         It is clear that
         $\gamma(\alpha)_{g,g}=1,\gamma(\alpha)_{g,h}=\gamma(\alpha)^{-1}_{h,g}.$
         \vskip 0.1in
         \noindent
         {\bf Lemma 3.1.14: }{\it Suppose that $gh=hg, gk=kg.$ Then
         \begin{description}
          \item[(1)] $\gamma(\delta \rho)_{g,h}=1$.

          \item[(2)] $\gamma(\alpha)_{g,
          hk}=\gamma(\alpha)_{g,h}\gamma(\alpha)_{g,k}.$
          \end{description}
          Hence, $L_g^{\alpha}=\gamma_{g,.}: C(g)\rightarrow U(1)$ is a representation of $C(g)$.}
          \vskip 0.1in
          \noindent
          {\bf Proof: } (1) is obvious. For (2),
          $$\begin{array}{lll}
           \gamma(\alpha)_{g,hk}&=&\alpha_{g,hk}\alpha_{hk,g}^{-1}\\
                                &=&\alpha_{g,hk}\alpha^{-1}_{gh,k}\alpha_{hg,k}
                                \alpha_{h,gk}^{-1}\alpha_{h,kg}\alpha^{-1}_{hk,g}\\
                                &=&\alpha_{g,h}\alpha^{-1}_{h,k}\alpha_{g,k}
                                    \alpha^{-1}_{h,g}\alpha_{h,k}\alpha^{-1}_{k,g}\\
                                &=&\gamma(\alpha)_{g,h}\gamma(\alpha)_{g,k}
           \end{array}
           $$
           \vskip 0.1in
           Recall the following classical definition.
           \vskip 0.1in
           \noindent
           {\bf Definition 3.1.15: }{\it $g$ is called
           $\alpha$-regular iff $L^{\alpha}_g$ is trivial.}
           \vskip 0.1in

            Next, we calculate discrete torsion for some groups.
    We first consider the case of a finite abelian group $G$. In this case
    $H^i(G, \Q)=0$ for $i\neq 0$. The exact sequence
    $$0\rightarrow \Z\rightarrow \C\rightarrow \C^*\rightarrow 1$$
    implies that $H^2(G,U(1))=H^2(G,\C^*)=H^3(G, \Z)$. By the universal coefficient
    theorem, $H^3(G,\Z)= H_2(G,\Z).$
    \vskip 0.1in
           \noindent
           {\bf Example 3.1.16 $G=\Z/n\times \Z/m$: } Note that $H^2(G, U(1))=H_2(G,\Z)=
    \Z/n\otimes \Z/m=Z_{gcd(n,m)}$. In this case, one can write down the phase of
    discrete torsion explicitly \cite{VW}. Let $\xi$ (resp. $\zeta$) be $n$ (resp. $m$) root
    of unity. Any element of $\Z/n\times \Z/m$ can be written as
    $(\xi^a, \zeta^b)$. Let $p=gcd(n,m)$. The phase of a discrete torsion can be
    written as
    $$\gamma_{(\xi^a, \zeta^b),(\xi^{a'},
    \zeta^{b'})}=\omega_p^{m(ab'-ba')}$$
    with $\omega_p=e^{2\pi i/p}, m=1, \cdots, p.$ There are
    $p$-phases for $p$-discrete torsions.  It is trivial
    to generalize this construction to an arbitrary finite abelian group.

    Suppose that $f: Y\rightarrow X$ is the orbifold universal
    cover and $G$ is the orbifold fundamental group which acts
    on $Y$ such that $X=Y/G$. Suppose $X_{(g)}$ is a
    sector (twisted or nontwisted) of $X$. For any $q\in X$, choose an orbifold chart $U_q/G_q$ satisfying
        Definition 3.1.8. A component of $f^{-1}(U_q/G_q)$ is of the form $U_q/\Gamma'$ for
    $\Gamma'\subset G_q$.
    It is clear that $G_q/\Gamma'$ is a subgroup of the orbifold fundamental
    group. Therefore, we obtain a group homomorphism
    $$\phi_q: G_q\rightarrow
\pi^{orb}_1(X).\leqno(3.1.18)$$
    It is easy to check that a different choice of component of $f^{-1}(U_q/G_q)$ or a
    different choice of $q\in X_{(g)}$ induces a
    homomorphism differing by a conjugation. Therefore, there is a unique
    map from the conjugacy classes of $G_q$ to the conjugacy classes of $\pi^{orb}_1(X)$.
    \vskip 0.1in
    \noindent
    {\bf Definition 3.1.17: }{\it We call $X_{(g)}$ a dormant sector if $\phi_p(g)=1$.}
    \vskip 0.1in
    If $X_{(g)}$ is a dormant sector, we define $L_{(g)}=1$. It will not
    receive any correction from discrete torsion. Non-dormant sectors are of
    the form $Y_g/C(g)$, where $Y_g\neq \emptyset$ is the fixed point locus of $1\neq g\in \pi^{orb}_1(X)$.
     $Y_{g}$ is a smooth suborbifold of $Y$.
    It is clear that $Y_{h^{-1}gh}$ is diffeomorphic to $Y_{g}$ by
    the action of $h$. By abusing the notation, we denote the twisted sector $Y_{g}/C(g)$ by
    $X_{(g)}$,
     where $C(g)$ is the centralizer of $g$.

    Let $\alpha$ be a  discrete
    torsion.  By Lemma 3.1.14(2),
      for each $g$, the phase
    $$L^{\alpha}_g: C(g)\rightarrow U(1)$$
    is a group homomorphism. We can use this group homomorphism to
    define a flat complex line-bundle
    $$L_{g}=Y_{g}\times_{L^{\alpha}_g}\C$$
    over $X_{(g)}$.
        \vskip 0.1in
    \noindent
    {\bf Lemma 3.1.18: }{\it
    \begin{description}
    \item[(1)] $L_{tgt^{-1}}$ is isomorphic to $L_{g}$ by the map
            $$t\times Id: Y_{g}\times \C\rightarrow Y_{tgt^{-1}}\times
                \C.\leqno(4.3)$$
                Hence, we can denote $L_{g}$ by $L_{(g)}$.
    \item[(2)] $L_{(g)}^{-1}=L_{(g^{-1})}.$
    \item[(3)] When we restrict to $X_{(g_1,\cdots,
    g_k)}=Y_{g_1}\cap\cdots \cap
                Y_{g_k}/C(g_1,\cdots, g_k)$,
            $L_{(g_1,\cdots, g_k)}=L_{(g_1)}\cdots L_{(g_k)}$, where
            $L_{(g_1,\cdots,  g_k)}=Y_{g_1}\cap\cdots
            \cap Y_{g_k}\times_{\gamma_{g_1\cdots g_k}} \C.$
    \end{description}}
    \vskip 0.1in
    \noindent
    {\bf Proof: } Recall that there is an isomorphism
    $$t_{\#}: C(g)\rightarrow C(tgt^{-1})$$
    given by $t_{\#}(h)=tht^{-1}$. The map
    $$t:Y_{g}\rightarrow X_{tgt^{-1}}$$ is $t_{\#}$-equivariant. By Lemma 3.8,
    $\gamma_{tgt^{-1}}(tht^{-1})=\gamma_{g}(h)$ for $h\in C(g)$. Then,
    $$(t\times Id)(hx, \gamma(h)(v))=(thx, \gamma_{g}(h)(v))=(tht^{-1}tx,
    \gamma_{tgt^{-1}}(tht^{-1})(v)).\leqno(3.1.19)$$
    Then we take the quotient by $C(g), C(tgt^{-1})$ respectively to get
    an isomorphism between $L_{g}, L_{tgt^{-1}}$.
    (2) and (3) follow from the fact that for any $h\in C(g_1,\cdots,g_k)$,
    $$\gamma(\alpha)_{g_1\cdots g_k, h}=\gamma(\alpha)^{-1}_{h,g_1\cdots g_k}
    =\gamma(\alpha)^{-1}_{h,g_1}\cdots \gamma(\alpha)^{-1}_{h, g_k}=
        \gamma(\alpha)_{g_1, h}\cdots\gamma(\alpha)_{g_k, h}.\leqno(3.1.20)$$

    \vskip 0.1in
    \noindent
    {\bf Theorem 3.1.19: }{\it $\L_{\alpha}=\{L_{(g)}\}_{(g)\in T_1}$ is an inner local system of $X$.}
    \vskip 0.1in
    \noindent
    {\bf Proof: } Property (1) is obvious. The property (2) follows from Lemma 3.1.17.
    Let's prove property (3).
         Consider the
    image $\g'=(g'_1, g'_2, g'_3)$ of $\g$ in $\pi^{orb}_1(X)$ under the homomorphism (3.1.18).
    Then we still have $g'_1g'_2g'_3=1$. There are three
    possibilities: (i) $g'_1=g'_2=g'_3=1$ and there is nothing to prove in
    this case; (ii) $g'_3=1, g'_2=(g')^{-1}_1$ is nontrivial; (iii)
    $g'_1,g'_2, g'_3$ are all nontrivial.

    For the second case, let $g=g'_1$. We have the following factorization
    $$e_1\times e_2\times e_3: X_{(\g)}\rightarrow X_{(g_1,g_2)}\times X_{(g_3)}
    \rightarrow X_{(g_1)}\times X_{(g_2)}\times X_{(g_3)}.$$
    However, $X_{(g_1, g_2)}=Y_{g}\cap Y_{(g^{-1})}/C(g,g^{-1})=Y_{g}/C(g).$
     Moreover, over $X_{g_1,g_2}$
    $$e^*_1 L_{(g)}e^*_2 L_{(g^{-1})}=L_{(g)}I^*L_{(g^{-1})}=1.\leqno(4.6)$$

    In the third case,
    $X_{(\g)}=Y_{g_1}\cap Y_{g_2}\cap Y_{g_3}/C(g_1, g_2, g_3)$. The
    proof follows from Lemma 3.1.18 (3).
    $\Box$
    \vskip 0.1in
    \noindent
    {\bf Remark 3.1.20: }{\it The current definition of discrete
    torsion is unsatisfactory because of the existence of dormant
    sectors (see example 3.5.4 as well). It would be desirable to find
    a better definition of discrete torsion where all the sectors
    will receive corrections.}
    \vskip 0.1in

\subsection{Degree shifting and orbifold cohomology group}
For the rest of the paper, we will assume that $X$ is an almost
complex orbifold with an almost complex structure $J$. Recall that
an almost complex structure $J$ on $X$ is a smooth section of the
orbifold bundle $End(TX)$ such that $J^2=-Id$. In this case,
multi-sectors $\widetilde{\Sigma_k X}$ naturally inherit an almost
complex structure. Moreover, both the evaluation map $e_{i_1,
\cdots, i_l}$ and $I$ are naturally pseudo-holomorphic, i.e., its
differential commutes with the almost complex structures on
$\widetilde{\Sigma_k X}$.

An important feature of orbifold cohomology groups is degree
shifting, which we shall explain now. Let $p\in X$ be a singular
point of $X$. The almost complex structure on $X$ gives rise to an
effective representation $\rho_p: G_p \rightarrow GL(n,\C)$ (here
$n=\dim_\C X$). For any $g\in G_p$, we write $\rho_p(g)$ as a
diagonal matrix $$ diag(e^{2\pi i m_{1,g}/m_g}, \cdots, e^{2\pi i
m_{n,g}/m_g}), $$ where $m_g$ is the order of $g$ in $G_p$, and
$0\leq m_{i,g} <m_g$. This matrix depends only on the conjugacy
class $(g)_{G_p}$ of $g$ in $G_p$. We define a function
$\iota:\widetilde{X}\rightarrow \Q$ by $$
\iota(p,(g)_{G_p})=\sum_{i=1}^n \frac{m_{i,g}}{m_g}. $$ It is
straightforward to show the following

\vspace{3mm}

\noindent{\bf Lemma 3.2.1: }{\it The function
$\iota:\widetilde{X}\rightarrow \Q$ is locally constant. We will
denote it by $\iota_{(g)}$. The function $\iota_{(g)}$ satisfies
the following conditions:

\begin{itemize}
\item $\iota_{(g)}$ is integral if and only if $\rho_p(g)\in SL(n,\C)$.
\item
$$ \iota_{(g)}+\iota_{(g^{-1})}=rank(\rho_p(g)-I), \leqno(3.2.1)$$
which is the ``complex codimension'' $\dim_\C X-\dim_\C
X_{(g)}=n-\dim_\C X_{(g)}$ of $X_{(g)}$ in $X$. As a consequence,
$\iota_{(g)}+\dim_\C X_{(g)}<n$.
\end{itemize}}

\vskip 0.1in

\noindent{\bf Definition 3.2.2: }{\it $\iota_{(g)}$ is called the
degree shifting number.}

\vskip 0.1in

In the definition of orbifold cohomology groups, we will shift up
the degree of cohomology classes of $X_{(g)}$ by $2\iota_{(g)}$.
The reason for such a degree shifting will become clear after we
discuss the dimension of the moduli space of ghost maps (see
Proposition 3.4.4).

    There are two important classes of orbifolds. $X$ is called an
    {\em $SL$-orbifold} if
    $\rho_p(g)\in SL(n, \C)$. $X$ is called an {\em $Sp$-orbifold} if
    $\rho_p(g)\in Sp(n,\C)$. In particular, a Calabi-Yau orbifold is a
    $SL$-orbifold. An holomorphic symplectic orbifold or
    hyperkahler orbifold is an $Sp$-orbifold.

    By the Lemma 3.2.1, $\iota_{(g)}$ is integral if  $X$ is a
    $SL$-orbifold.

We observe that although the almost complex structure $J$ is
involved in the definition of degree shifting numbers
$\iota_{(g)}$, they do not depend on $J$ because locally the
parameter space of almost complex structures, which is the coset
$SO(2n,\R)/U(n,\C)$, is connected.

\vskip 0.1in

\noindent{\bf Definition 3.2.3: }{\it Let $\L$ be an inner local
system. We define the orbifold cohomology groups $H^d_{orb
}(X;\L)$ and compactly supported orbifold cohomology group
$H^d_{orb, c}(X, \L)$ of $X$ by $$ H^d_{orb}(X;\L)=\oplus_{(g)\in
T} H^{d-2\iota_{(g)}}(X_{(g)};\L_{(g)}),
H^d_{orb,c}(X;\L)=\oplus_{(g)\in T}
H^{d-2\iota_{(g)}}_c(X_{(g)};\L_{(g)})\leqno(3.2.2) $$ and
orbifold Betti numbers $b^d_{orb, \L}=\sum_{(g)}\dim
H^{d-2\iota_{(g)}}
    (X_{(g)},\L_{(g)})$. If $\L=\L_{\alpha}$ for some discrete torsion
    $\alpha$, we define $H^*_{orb, \alpha}(X, \C)=H^*_{orb}(X,
    \L_{\alpha}), H^*_{orb,c, \alpha}(X, \C)=H^*_{orb,c}(X,
    \L_{\alpha}) $.}

\vskip 0.1in

Note that, in general, orbifold cohomology groups are rationally
graded. Traditionally, $H^*_{orb}(X, \L)$ for $\L=1$ is called
ordinary orbifold cohomology. Other cases are called twisted
orbifold cohomology.

\vspace{3mm}

Suppose $X$ is a  complex orbifold with an integrable complex
structure $J$. Then each twisted sector $X_{(g)}$ is also a
complex orbifold with the induced complex structure. We consider
the \v{C}ech cohomology groups on $X$ and on each $X_{(g)}$ with
coefficients in the sheaves of holomorphic forms (in the orbifold
sense). These \v{C}ech cohomology groups are identified with the
Dolbeault cohomology groups of $(p,q)$-forms (in the orbifold
sense). When $X$ is closed, the harmonic theory  can be applied to
show that these groups are finite dimensional, and there is a
Kodaira-Serre duality between them. When $X$ is a closed K\"ahler
orbifold (so is each $X_{(g)}$), these groups are then related to
the singular cohomology groups of $X$ and $X_{(g)}$ as in the
smooth case, and the Hodge decomposition theorem holds for these
cohomology groups.

\vspace{3mm}

\noindent{\bf Definition 3.2.4:} {\it Let $X$ be a closed complex
orbifold. We define, for $0\leq p,q\leq \dim_{\C}X$, orbifold
Dolbeault cohomology groups $$ H^{p,q}_{orb}(X;\L)=\oplus_{(g)}
H^{p-\iota_{(g)},q-\iota_{(g)}} (X_{(g)};\L_{(g)}),
H^{p,q}_{orb,c}(X;\L)=\oplus_{(g)}
H^{p-\iota_{(g)},q-\iota_{(g)}}_c (X_{(g)};\L_{(g)}).
\leqno(3.2.3) $$ We define orbifold Hodge numbers by
$h^{p,q}_{orb}(X, \L)=\dim H^{p,q}_{orb}(X;\L_{(g)})$.}

\vspace{3mm}
    \noindent
    {\bf Remark 3.3.3: }{\it In the case of global quotient
    $X=Y/G$, suppose that $\alpha\in H^2(G, S^1)$. $L^{\alpha}_g$
    induces a twisted action of $C(g)$ on the fixed point set
    $H^*(Y_g, \C)$ by $h\circ \beta=L^{\alpha}_g (h) h^*\beta$.
    Let $H^*(Y_g, \C)^{C^{\alpha}(g)}$ be the invariant subspace under
    the twisted action. It is easy to observe that
    $$H^d_{orb, \alpha}(X. \C)=\oplus_{(g)} H^{d-2\iota_{(g)}}(Y_g,
    \C)^{C^{\alpha}(g)}$$
    and
    $$ H^{p,q}_{orb, \alpha}(X;\C)=\oplus_{(g)}
H^{p-\iota_{(g)},q-\iota_{(g)}} (Y_g, \C)^{C^{\alpha}(g)}.$$}
\vskip 0.1in

\subsection{Poincar\'e duality}
    For simplicity, we assume that the orbifold under consideration is closed.

Recall that there is a natural $C^\infty$ map
$I:X_{(g)}\rightarrow X_{(g^{-1})}$ defined by $(p,(g))\rightarrow
(p,(g^{-1}))$, which is an automorphism of $\widetilde{X}$ as an
orbifold, and $I^2=Id$ (Remark 3.1.4).

\vspace{3mm}

\noindent{\bf Proposition 3.3.1: }(Poincar\'{e} duality){\it

For any $0\leq d\leq 2n$, the pairing
 $$ <\,\ >_{orb}:
H^d_{orb}(X;\L)\otimes H^{2n-d}_{orb,c}(X;\L)\rightarrow\C $$
defined by the direct sum of
    $$<\,\ >_{orb}^{(g)}:
H^{d-2\iota_{(g)}}(X_{(g)};\L_{(g)})\otimes
H^{2n-d-2\iota_{(g^{-1})}}_c(X_{(g^{-1})};\L_{(g^{-1})})\rightarrow\C
$$ where
    $$<\alpha, \beta>^{(g)}_{orb} =\int_{X_{(g)}} \alpha\wedge
I^*(\beta)\leqno(3.3.3)$$
    for $\alpha\in H^{d-2\iota_{(g)}}(X_{(g)}, \L_{(g)}), \beta\in
    H^{2n-d-2\iota_{(g^{-1})}}(X_{(g^{-1})}, \L_{(g^{-1})})$
    is nondegenerate. }

\vskip 0.1in

By Lemma 3.1.18, $I^*L_{(g^{-1})}=L^{-1}_{(g)}$. Hence, the
definition makes sense. Moreover,  $<\,\
>_{orb}$ is just  the ordinary Poincar\'{e} pairing when
restricted to the nontwisted sector. \vskip 0.1in

\noindent{\bf Proof:}
    By (3.2.1), we have
    $$
     2n-d-2\iota_{(g^{-1})}=\dim X_{(g)}-d-2\iota_{(g)}.\leqno(3.3.4)
    $$
    Furthermore, $I|_{X_{(g)}}: X_{(g)}\rightarrow X_{(g^{-1})}$ is
    a homeomorphism and $I^*L_{(g^{-1})}=L^{-1}_{(g)}$.
    Under this homeomorphism, $<\,\ >^{(g)}_{orb}$
    is isomorphic to the ordinary Poincar\'{e} pairing with coefficients on $X_{(g)}$.
    Hence $<\,\ >_{orb}$ is nondegenerate.
    \hfill $\Box$

\vspace{3mm}

For the case of orbifold Dolbeault cohomology, the following
proposition is straightforward.

\vspace{3mm}

\noindent{\bf Proposition 3.3.2:} {\it Let $X$ be an
$n$-dimensional  complex orbifold. There is a Kodaira-Serre
duality pairing $$ <\,\ >_{orb}:H^{p,q}_{orb}(X;\L)\times
H^{n-p,n-q}_{orb,c}(X;\L)\rightarrow \C $$ similarly defined as in
the previous proposition. When $X$ is closed and K\"{a}hler, the
following is true:
\begin{itemize}
\item $H^r_{orb}(X;\L)=\oplus_{r=p+q}H^{p,q}_{orb}(X;\L)$
\item $H^{p,q}_{orb}(X;\L)=\overline{H^{q,p}_{orb}(X;\L)}$,
\end{itemize}
and the two pairings (Poincar\'{e} and Kodaira-Serre) coincide.}
\vskip 0.1in

\subsection{Orbifold cup product}
    Our definition of orbifold cup product is motivated by the
construction of quantum product. For this approach, we  have to
construct a Gromov-Witten type invariant of genus zero,  homology
class zero with three marked points. It involves an analysis of
the moduli space of constant (ghost) good maps from the orbifold
sphere with three marked points. The construction is lengthy.
However, we can skip over the construction of the moduli space and
write down the definition explicitly. If the reader wishes to
understand the geometric origin of our definition and key
properties such as associativity, we encourage the reader to read
\cite{CR1}.

    A key ingredient is  orbifold Riemann surface. Every
closed orbifold of dimension 2 is complex, with underlying
topological space a closed Riemann surface. More precisely, a
closed 2-dimensional orbifold consists of the following data: a
closed Riemann surface $\Sigma$ with complex structure $j$, a
finite subset of distinct points $\z= (z_1,\cdots,z_k)$ on
$\Sigma$, each with  multiplicity $m_i\geq 2$ (let
$\m=(m_1,\cdots,m_k)$), such that the orbifold structure at $z_i$
is given by the ramified covering $z\rightarrow z^{m_i}$. We will
also call a closed 2-dimensional orbifold a {\it complex
orbicurve} when the underlying complex analytic structure is
emphasized.

We observe the following well-known fact (see e.g. \cite{SC}).

\vspace{3mm}

\noindent{\bf Proposition 3.4.1: }{\it Let $(\Sigma,\z,\m)$ be a
complex orbicurve, where $\z=(z_1,\cdots,z_k)$ and
$\m=(m_1,\cdots,m_k)$, such that either $g_\Sigma\geq 1$, or
$g_\Sigma=0$ with $k\geq 3$ or $k=2$ and $m_1=m_2$. Then there is
a closed Riemann surface $\widetilde{\Sigma}$ with a finite group
$G$ acting on $\widetilde{\Sigma}$ holomorphically, and a
holomorphic mapping $\pi:\widetilde{\Sigma}\rightarrow \Sigma$,
such that $(\widetilde{\Sigma},G,\pi)$ is a uniformizing system of
$(\Sigma,\z,\m)$.}

\vspace{3mm}

The construction of cup product follows the procedure to define
   quantum product. First, we need to define a 3-point function.
   In our case, the 3-point function is an integral over $X_{(\g)}$
   for $(\g)\in T^o_3$. In order to write down the form we
   integrate, we need to construct an obstruction bundle and its
   Euler form over $X_{(\g)}$.

   Consider the pull-back tangent bundle $e^*TX$ over $X_{(\g)}$.
   Let $x\in X_{(\g)}$ be a generic point and its local group in
   $X$ is $G'$. It obviously contains three elements $g_1, g_2, g_3$ with the relation
   $g_1g_2g_3=1, g^{k_i}_i=1$, where $k_i$ is the order of $g_i$.
   Let $G$ be the subgroup of $G'$ generated by $g_1, g_2, g_3$.
   Clearly, $G$ acts on $e^*TX$ while fixing $X_{(\g)}$

   Consider an orbifold Riemann sphere with three orbifold points
   $(S^2, (x_1, x_2, x_3), (k_1, k_2, k_2))$. Without any confusion, we simply denote it by
   $S^2$.
   Recall Example 3.1.12
   $$\pi^{orb}_1(S^2)=\{\lambda_1, \lambda_2, \lambda_3;
   \lambda^{k_i}_i=1, \lambda_1\lambda_2\lambda_3=1\},$$
   where $\lambda_i$ is represented by a loop around the marked point
   $x_i$. There is an obvious surjective homomorphism
   $$\pi: \pi^{orb}_1(S^2)\rightarrow G.\leqno(3.1.4)$$
   $ker \pi$ is a subgroup of finite index. Suppose that
   $\tilde{\Sigma}$ is the orbifold universal cover of $S^2$.
   By  Proposition
   3.4.1, it is uniformized by a closed Riemann surface $\Sigma'$.
   Hence, $\tilde{\Sigma}$ is the universal cover of $\Sigma'$ and
   is smooth. Let $\Sigma=\tilde{\Sigma}/ker \pi$. $\Sigma$ is
   compact and $S^2=\Sigma/G$. Since $G$ contains the relation
   $g^{k_i}_i=1$, $\Sigma$ is smooth.
    $G$ acts on
   $H^1(\Sigma).$ Let $e: X_{(\g)}\rightarrow X$ be the evaluation
   map.  Therefore, we can assume
   that $G$ acts on both $H^1(\Sigma)$ and $e^*TX$. We view $H^1(\Sigma)$ as
   a trivial bundle over $X_{(\g)}$. The obstruction bundle $E_{(\g)}$ we want is
   is the invariant part of $H^1(\Sigma)\otimes e^*TX$, i.e.,
   $E_{(\g)}=(H^1(\Sigma)\otimes e^*TX)^G$. Since we do not assume that $X$ is compact, $X_{(\g)}$
   could be a non-compact orbifold in general. The Euler class of $E_{(\g)}$ depends on a choice of connection
   on $E_{(\g)}$.  Let $e_A(E_{(\g)})$ be the Euler form
   computed from the connection $A$ by Chern-Weil theory. It is clear
   that $e_A(E_{(\g)}), e_{A'}(E_{(\g)})$ differ by an exact form if $A'$ is
   another connection on $E_{(\g)}$.

    Now, we are ready to define our 3-point function. Suppose that
    $\alpha \in H^{d_1}_{orb}(X, \C), \beta\in H^{d_2}_{orb}(X, \C),
    \gamma\in H^{*}_{orb, c}(X_{(g_3}, \C)$.
    \vskip 0.1in
    \noindent
    {\bf Definition 3.1.3: }{\it We define the 3-point function
    $$<\alpha, \beta, \gamma>_{orb}=\sum_{(\g)\in T^0_3}\int_{X_{(\g)}}
    e^*_1\alpha\wedge e^*_2\beta\wedge e^*_3\gamma \wedge
    e_A(E_{(\g)}).\leqno(3.1.5)$$}
    \vskip 0.1in
    Note that $e^*_3\gamma$ is compact supported. Therefore, the
    integral is finite. Moreover, if we choose a different $A'$,
    then $e_A(E_{(\g)}), e_{(A'}(E_{(\g)})$ differ by an exact form. Hence, the
    integral is independent of the choice of $A$.

    \vskip 0.1in
    \noindent
    {\bf Definition 3.1.4: }{ \it We define the orbifold cup product
    by the relation
    $$<\alpha\cup_{orb} \beta, \gamma>_{orb}=<\alpha, \beta,
    \gamma>_{orb}.\leqno(3.1.6)$$
    }
    \vskip 0.1in
    If $\alpha, \beta$ are compacted supported orbifold cohomology classes,
    we can define $\alpha\cup_{orb} \beta\in H^*_{orb, c}(X, \L)$ in the same fashion.
    Suppose that $\alpha\in H^*(X_{(g_1)}, \C), \beta\in
    H^*(X_{(g_2)}, \C)$. $\alpha\cup_{orb}\beta\in H^*_{orb}(X,
    \C)=\oplus_{(g)\in T_1} H^*(X_{(g)}, \C)$. Therefore, we
    should be able to decompose $\alpha\cup_{orb}\beta$ as a sum
    of its components in $H^*(X_{(g)}, \C)$. Such a decomposition
    would be very useful in computation.

    Note that when $g_1g_2g_3=1$ the conjugacy class $(g_1, g_2,
    g_3)$ is uniquely determined by the conjugacy class of the pair
    $(g_1, g_2)$. We can use it to obtain the following
    \vskip 0.1in
    \noindent
    {\bf Decomposition Lemma 3.1.5: }{\it
    $$\alpha\cup_{orb}\beta =\sum_{(h_1, h_2)\in T_2, h_i\in
    (g_i)} (\alpha\cup_{orb}\beta)_{(h_1,h_2)},\leqno(3.1.7)$$
    where $(\alpha\cup_{orb}\beta)_{(h_1,h_2)}\in
    H^*(X_{(h_1h_2)}, \C)$ is defined by the relation
    $$<(\alpha\cup_{orb}\beta)_{(h_1,h_2)}, \gamma>_{orb}
    =\int_{X_{(h_1,h_2)}}e^*_1\alpha\wedge e^*_2\beta\wedge e^*_3\gamma \wedge
    e_A(E_{(\g)}).\leqno(3.1.8)$$
    for $\gamma\in H^*_c(X_{((h_1h_2)^{-1})}, \C)$.}
    \vskip 0.1in
    \noindent
    {\bf Remark: }{\it In the case of global quotient $X=Y/G$,
    $X_{(h_1, h_2)}=Y_{h_1}\cap Y_{h_2}/ C(h_1)\cap C(h_2)$.}
    \vspace{3mm}

\noindent{\bf Theorem 3.1.6: } {\it Let $X$ be an almost complex
orbifold with almost complex structure $J$ and $\dim_\C X=n$. The
cup product defined above preserves the orbifold degree
$\cup_{orb}: H_{orb}^p(X;\F)\otimes H_{orb}^q(X;\F) \rightarrow
H_{orb}^{p+q}(X;\F)$ for any $0\leq p,q\leq 2n$ such that $p+q\leq
2n$, and has the following properties:
\begin{enumerate}
\item The total orbifold cohomology group $H^\ast_{orb}(X;\F)=
\oplus_{0\leq d\leq 2n}H^d_{orb}(X;\F)$ is a ring with unit
$e_X^0\in H^0(X;\F)$ under $\cup_{orb}$, where $e_X^0$ is the
Poincar\'{e} dual to the fundamental class $[X]$.

\item The cup product $\cup_{orb}$ is invariant under deformation of
$J$.
\item When $X$ is of integral degree shifting numbers, the total orbifold
cohomology group $H_{orb}^\ast(X;\F)$ is integrally graded, and we
have supercommutativity $$
\alpha_1\cup_{orb}\alpha_2=(-1)^{\deg\alpha_1\cdot\deg\alpha_2}\alpha_2
\cup_{orb}\alpha_1. $$
\item Restricted to the nontwisted sectors, i.e., the ordinary
cohomology $H^\ast(X;\F)$, the cup product $\cup_{orb}$ equals the
ordinary cup product on $X$.
\item $\cup_{orb}$ is associative.
\end{enumerate}}

\vspace{3mm}

\vspace{3mm}

Now we define the cup product $\cup_{orb}$ on the total orbifold
Dolbeault cohomology group of $X$ when $X$ is a  complex orbifold.
We observe that in this case all the objects we have been dealing
with are holomorphic, i.e., $\widetilde{\Sigma_k X}$ is a complex
orbifold, $pr:E_{(\g)}\rightarrow X_{(\g)}$ is holomorphic
orbifold bundle, and the evaluation map are holomorphic.

\vspace{3mm}

\noindent{\bf Definition 4.1.7: }{\it For any $\alpha_1\in
H^{p,q}_{orb}(X;\C)$, $\alpha_2\in
H_{orb}^{p^\prime,q^\prime}(X;\C)$, we define the 3-point function
and orbifold cup product in the same fashion as Definition 3.1.3,
3.1.4.} \hfill $\Box$

\vskip 0.1in

Note that since the top Chern class of a holomorphic orbifold
bundle can be represented by a closed $(r,r)$-form, where $r$ is
the rank, it follows that $\alpha_1\cup_{orb}\alpha_2$ lies in
$H_{orb}^{p+p^\prime,q+q^\prime}(X;\C)$.

\vspace{3mm}

The following theorem can be similarly proved.

\vspace{3mm}

\noindent{\bf Theorem 4.1.8: }{\it Let $X$ be an $n$-dimensional
closed complex orbifold with complex structure $J$. The orbifold
cup product $$ \cup_{orb}: H_{orb}^{p,q}(X;\C)\otimes
H_{orb}^{p^\prime,q^\prime}(X;\C) \rightarrow
H_{orb}^{p+p^\prime,q+q^\prime}(X;\C) $$ defined above has the
following properties:
\begin{enumerate}
\item The total orbifold Dolbeault cohomology group is a ring with unit
$e_X^0\in H_{orb}^{0,0}(X;\C)$ under $\cup_{orb}$, where $e_X^0$
is the class represented by the equaling-one constant function on
$X$.

\item The cup product $\cup_{orb}$ is invariant under
deformation of $J$.
\item When $X$ is of integral degree shifting numbers, the total orbifold
Dolbeault cohomology group of $X$ is integrally graded, and we
have supercommutativity $$
\alpha_1\cup_{orb}\alpha_2=(-1)^{\deg\alpha_1\cdot\deg\alpha_2}\alpha_2
\cup_{orb}\alpha_1. $$
\item Restricted to the nontwisted sectors, i.e., the ordinary Dolbeault
cohomology $H^{\ast,\ast}(X;\C)$, the cup product $\cup_{orb}$
equals the ordinary wedge product on $X$.
\item When $X$ is K\"{a}hler and closed, the cup product $\cup_{orb}$ coincides with
the orbifold cup product over the orbifold cohomology groups
$H^\ast_{orb}(X;\F)$ under the relation $$
H^r_{orb}(X;\F)=\oplus_{p+q=r}H^{p,q}_{orb}(X;\C). $$
\end{enumerate}}
    \vskip 0.1in
    \noindent
    The most difficult part of proof is associativity. We refer the proof to \cite{CR1}.

    \vskip 0.1in
    \noindent
    {\bf Remark 3.4.9: }{\it

    In many way, the current definition of orbifold cohomology is
    less than
    satisfactory. It is a very interesting question whether one
    can represent a cohomology class from a twisted sector by a
    differential form on $X$ with certain singularities along
    singular strata. Recall that a crepant resolution is a map
    $F: Y\rightarrow X$ such that $F^*K_X=K_Y$, where $X$ is a
    complex orbifold and $Y$ is a smooth complex manifold. An
    extremely interesting question is to study the relation
    between the orbifold cohomology of $X$ and the ordinary cohomology of
    $Y$ (see section 6). The main difficulty is to pull back the class from
    the twisted sector. If we can establish a theory to use
    differential forms on $X$ to represent the class from a twisted
    sector, it would be very useful to understand its relation to
    crepant resolution.

    Note that orbifold cohomology is a cohomology theory only. It would be
    desirable to build a homology theory based on suborbifolds and their intersection
    theory. Hence, we can realize orbifold Poincar\'{e} duality and orbifold cup
    product geometrically. It seems to the author that good maps should
    play a critical role in such a homology theory.}

\subsection{Examples}
    So far, only a few examples of global quotients have been computed by
physicists \cite{VW}  \cite{D}. However, orbifold cohomology is
very much calculable, as we will demonstrate in examples. Here we
compute
    several  examples. The first two are  local examples, where the reader should have some
    general ideas about orbifold cohomology. The third and fourth have nontrivial discrete torsion.
    One is a global quotient and
    another one is a non-global quotient. The fourth example has the phenomenon that  most
    of the twisted sectors are dormant sectors. The last one is Joyce's \cite{JO} example, where there is
    no nontrivial discrete torsion. However, there are nontrivial inner local systems. We will compute
    the twisted orbifold cohomology given by  nontrivial inner local systems to match
    Joyce's desingularizations.
    \vskip 0.1in
    \noindent
    {\bf Example 3.5.1: } The easiest example is probably a point
    with a trivial group action of $G$. In this case, the orbifold cohomology is
    generated by conjugacy classes of elements of $G$. All the
    degree shifting numbers are zero. Only the Poincar\'{e} paring and cup
    products are interesting. We observe that
    $X_{(g_1,g_2 (g_1g_2)^{-1})}$ is a point with a trivial group action of $C(g_1)\cap C(g_2)$. Hence,
    the GW-invariants
    $$\int_{X_{(g_1,g_2 (g_1g_2)^{-1})}} 1=\frac{1}{|C(g_1)\cap C(g_2)|}.$$

    Let $x_{(g)}$ be the generators of orbifold cohomology group corresponding to one on sector $X_{(g)}$. Using decomposition lemma 4.5.1, the cup product
    $$x_{(g_1)}\cup x_{(g_2)}=\sum_{(h_1, h_2), h_1\in (g_1), h_2\in (g_2)}d_{(h_1,h_2)}x_{(h_1h_2)},$$
    where $(h_1,h_2)$ is the conjugacy class of pair $h_1, h_2$.
    We can compute coefficient $d_{(h_1,h_2)}$ by the relation
    $$\frac{1}{|C(g_1)\cap C(g_2)|}=\int_{X_{(h_1,h_2
    (h_1h_2)^{-1})}}
    1=d_{(h_1,h_2)}\int_{X_{(h_1h_2)}}1=\frac{d_{(h_1,h_2)}}{|C(h_1h_2)|}.$$
    Hence, $d_{(h_1,h_2)}=\frac{|C(h_1h_2)|}{|C(h_1)\cap C(h_2)|}$.
    Recall that the center $Z(\C[G])$ of group algebra $\C[G]$ is generated
    by $\sum_{h\in (g)} h$ and hence can be identified with conjugacy classes itself. In this sense,
    we say the center is generated by the set of conjugacy classes. $Z(\C[G])$ has
    a natural ring structure inherited from $\C[G]$.  Note that if we enumerate all elements in
    conjugacy class of pair $(h_1, h_2)$. Its product run through the conjugacy class $(h_1h_2)$
    precisely $\frac{|C(h_1h_2)|}{|C(h_1)\cap C(h_2)|}$ times.
    Therefore, the
    orbifold cup product  is the same as product of $Z(\C[G])$.

     Suppose
    that $\alpha\in H^2(G, U(1))$ is a discrete torsion. It is
    clear that twisted orbifold cohomology is generated by
    conjugacy classes of   $\alpha$-regular elements. On the
    another hand, the center of the twisted group algebra $\C_{\alpha}[G]$
    is also generated by conjugacy classes of $\alpha$-regular
    elements. Indeed, they have the same ring structures.

    \vskip 0.1in
    \noindent
    {\bf Example 3.5.2: } Suppose that $G\subset GL(n, \C)$ is a
    finite subgroup. Then
    $\C^n/G$ is an orbifold. Suppose that $\alpha\in H^2(G, U(1))$
    is a discrete torsion. For any $g\in G$, the fixed point set
    $X_g$ is a vector subspace and  $X_{(g)}=X_g/C(g)$. By
    definition, $L_{(g)}=X_g\times_{L^{\alpha}_g} \C$.
    Therefore, $H^*(X_{(g)}, L_{(g)})$ is the  subspace of
    $H^*(X_g, \C)$ invariant under the twisted action of $C(g)$
    $$h\circ \beta=L^{\alpha}_g(h) h^*\beta\leqno(3.5.2)$$
    for any $h\in C(g), \beta\in H^*(X_g, \C)$. However,
    $H^i(X_g, \C)=0$ for $i\geq 1$. Moreover, if
    $L^{\alpha}_g$ is nontrivial, $H^0(X_g, L_{(g)})=0$.
    Therefore,
    $H^{p, q}_{orb}=0$ for $p\neq q$ and
    $H^{p,p}_{orb}$ is
    a vector space generated by conjugacy classes of
    $\alpha$-regular elements $g$ with $\iota_{(g)}=p$.
    Therefore, we have a natural decomposition
    $$Z[\C_{\alpha}[G])=\sum_p H_p,\leqno(3.5.3)$$
    where $H_p$ is generated by conjugacy classes of
    $\alpha$-regular elements $g$ with $\iota_{(g)}=p$.
    The ring structure is also easy to describe. Let $x_{(g)}$ be
    the
    generator corresponding to the zero cohomology class of the twisted
    sector $X_{(g)}$ such that $g$ is $\alpha$-regular. We would
    like to get a formula for $x_{(g_1)}\cup x_{(g_2)}$. As we showed before,
    the multiplication of conjugacy classes can be described in terms
    of the
    center of the twisted group algebra $Z(\C_{\alpha}[G])$. But we have further
    restrictions in this case. Let's first describe the moduli space
    $X_{(h_1,h_2,(h_1h_2)^{-1})}$ and its corresponding GW-invariants. It is
    clear that
    $$X_{(h_1,h_2,(h_1h_2)^{-1})}=X_{h_1}\cap X_{h_2}/C(h_1)\cap C(h_2).$$
    To have the nonzero invariant, we require that
    $$\iota_{(h_1h_2)}=\iota_{(h_1)}+\iota_{(h_2)}.\leqno(3.5.4)$$
     Then we need to compute
   $$\int_{X_{h_1}\cap X_{h_2}/C(h_1,
   h_2)}e^*_3(vol_c(X_{h_1h_2}/C(h_1 h_2)))\wedge e(E),\leqno(3.5.6.1)$$
   where $vol_c(X_{h_1h_2}/C(h_1 h_2))$ is the compactly supported $C(h_1h_2)$-invariant
   top form with volume one
   over $X_{h_1h_2}$.
   However,
   $$X_{h_1}\cap X_{h_2}\subset X_{h_1h_2}$$
   is a submanifold. Therefore, the integral (3.5.5) is zero unless
   $$X_{h_1}\cap X_{h_2}=X_{h_1h_2}.\leqno(3.5.6.2)$$
   In this case, we call $(h_1,h_2)$ transverse. In this case,
   it is clear that the obstruction bundle is trivial. Suppose
   that
   the  integral is $d_{h_1, h_2}$. Let
   $$I_{g_1, g_2}=\{(h_1, h_2); h_i \in (g_i),
   \iota_{(h_1)}+\iota_{(h_2)}=\iota_{(h_1h_2)}, (h_1,
   h_2)-transverse \}.\leqno(3.5.6.3)$$
   Then,
   $$x_{(g_1)}\cup x_{(g_2)}=\sum_{(h_1,h_2)\in I_{g_1, g_2}}
   d_{h_1,h_2} x_{(h_1h_2)}.\leqno(3.5.6.4)$$
   A similar computation as previous example yields $d_{(h_1,h_2)}=
     \frac{|C(h_1h_2)|}{|C(h_1)\cap C(h_2)|}$.

   Next, we specialize to the symmetric group $S_n$. Recall that any element of the symmetric group
   can be decomposed into cycles $\gamma_1, \cdots, \gamma_m$ and its conjugacy class is uniquely
   determined by the cycle classes. Following  tradition, for each cycle $\gamma$
   of $k$ letters, we define its degree
   $deg(\gamma)=k-1$. In particular, the identity element has degree zero. Such a definition of
   degree is also invariant under the inclusion $S_n\subset S_{n+1}$. Then we define the degree of
   a conjugacy class as the sum of the degree of its cycles. It is easy to compute that the
   degree shifting number is the degree.  Moreover, we observe that
   $\iota_{(h_1h_2)}\leq \iota_{(g_1)}+\iota_{(g_2)}$ for $h_i\in (g_i)$. When equality holds,
   the pair $(h_1, h_2)$ is    automatically transverse.
   Hence,
   $$x_{(g_1)}\cup x_{(g_2)}=P_{\iota_{(g_1)}+\iota_{(g_2)}}(x_{(g_1)}x_{(g_2)}),\leqno(3.5.6.5)$$
   where $P_p: Z(\C_{\alpha}[G])\rightarrow H_p$ is the
   projection.
   The formula is precisely the formula appeared in Lehn-Sorger's
   calculation of the cohomology ring of the Hilbert scheme of points of
   $\C^2$. Therefore, by combining with their calculation, we have
   \vskip 0.1in
   \noindent
   {\bf Corollary 3.5.2a: }{\it $H^*(Hilb_k (\C^2), \C)$ and
   $H^*_{orb}(Sym_k (\C^2, \C)$ are isomorphic as rings.}
   \vskip 0.1in
   We would like to emphasize that our calculation works over  any finite quotient of affine space
   and is much more general than the symmetric
   group.

    \vskip 0.1in
    It is also easy to compute ring structure for the following
    examples. We leave it to readers.
    \vskip 0.1in
    \noindent
    {\bf Example 3.5.3 $T^4/\Z_2\times \Z_2$: } Here, $T^4=\C^2/\bigwedge$, where $\bigwedge$ is
    the lattice of integral points. Suppose that $g,h$ are generators of
    the first and the second factor of $\Z_2\times \Z_2$.
        The action of $\Z_2\times \Z_2$ on $T^4$ is defined as
    $$g(z_1, z_2)=(-z_1, z_2), h(z_1, z_2)=(z_1, -z_2).\leqno(3.5.6.6)$$
    The fixed point locus of $g$ is 4 copies of $T^2$. When we divide it by the remaining action
    generated by $h$, we obtain twisted sectors consisting of 4 copies of $S^2$. The
    degree shifting number for these twisted sectors is $\frac{1}{2}$. For the same reason,
    the fixed point locus of $h$ gives twisted sectors consisting of $4$ copies of
    $S^2$ with degree shifting number $\frac{1}{2}$. The fixed point locus of $gh$ is
    16 points, which are fixed by the whole group. The degree shifting number of the 16 points
    is $1$. An easy calculation shows that the nontwisted sector contributes one generator to degree
    0, 4 orbifold cohomology and two generators to degree 2 orbifold cohomology and no other.
    Using this information, we can compute the ordinary orbifold cohomology group
    $$b_{orb}^0=b_{orb}^4=1, b_{orb}^1=b_{orb}^3=8, b_{orb}^2=18.\leqno(3.5.7)$$

    By example 2.10, $H^2(\Z_2\times \Z_2, U(1))=\Z_2$. By Remark 2.2, the nontrivial generator of
    $H^2(\Z_2\times \Z_2, U(1))$ induces a discrete torsion $\alpha$. Next, we compute the
    twisted orbifold cohomology $H^*_{orb, \alpha}(T^4/\Z_2\times \Z_2, \C)$. Note that
    $\gamma (\alpha)_{gh,g}=\gamma (\alpha)_{gh,h}=-1$. Hence, the flat orbifold line bundles over
    the twisted sectors given by the 16 fixed points of $gh$ are nontrivial. Therefore, they
    contribute nothing to twisted orbifold cohomology.
    For  two-dimensional twisted sectors, let's consider a
    component of the fixed point locus of $g$. By the previous description, it is $T^2$. $h$ acts on
    $T^2$. Then the twisted sector $S^2=T^2/\{h\}$. We observe that the flat orbifold line bundle
    over $S^2$ is constructed as $L=T^2\times_{\gamma(\alpha)_g} \C$. Hence $H^*(S^2, L)$
    is isomorphic to the space of invariant cohomology of $T^2$ under the action    of $h$ twisted by $\gamma(\alpha)_g$ as $h(\beta)=\gamma(\alpha)_{g,h}h^*\beta$.
    By example 2.10, $\gamma(\alpha)_{g,h}=-1$. The invariant cohomology is $H^1(T^2, \C)$. Using the degree
    shifting number to shift up its degree, we obtain the twisted orbifold cohomology
    $$b_{orb,\alpha}^0=b_{orb,\alpha}^4=1, b_{orb,\alpha}^1=b_{orb,\alpha}^3=0,     b_{orb,\alpha}^2=18.\leqno(3.5.8)$$

    \vskip 0.1in
    \noindent
    {\bf Example 3.5.4 $WP(2, 2d_1)\times WP(2, 2d_2)$ ($d_1,d_2>1, (d_1,d_2)=1$): } Here, $WP(2, 2d)$ is the weight projective
    space of weighted $(2, 2d)$. $WP(2, 2d_1)\times WP(2, 2d_2)$ is not a global quotient unless
    $d_1=d_2=1$. In fact, its orbifold universal cover is $WP(1, d_1)\times WP(1, d_2)$ and
    $WP(2, 2d_1)\times WP(2, 2d_2)=WP(1, d_1)\times WP(1, d_2)/\Z_2\times \Z_2$. Hence, the orbifold
    fundamental group is $\Z_2\times \Z_2$. Therefore, there is
    a nontrivial discrete torsion $\alpha\in H^2(\Z_2\times \Z_2, U(1))$.

Next, we describe the
    twisted sectors. Suppose that $p=[0,1],q=[1,0]\in WP(1,d_1)$. We also
    use $p,q$ to denote its image in $WP(2,2d_1)$.  We use $p',q'$ to denote
    the corresponding points in $WP(1,d_2), WP(2,2d_2)$. $\{p\}\times WP(2,2d_2), \{p'\}\times
    WP(2,2d_1)$ give rise to two twisted sectors with degree shifting number $\frac{1}{2}$.
    $\{q\}\times WP(2,2d_2),\{q'\}\times WP(2,2d_1)$ give rise to $2d_1-1,2d_2-1$  twisted sectors
    with degree shifting numbers $\frac{i}{2d_1},\frac{j}{2d_2}$ for $1\leq i\leq
    2d_1-1,1\leq j\leq 2d_2-1$.
    $\{p\}\times\{p'\}$ gives rise to a twisted sector with degree shifting number $1$.
    $\{p\}\times \{q'\}$ gives rise to $2d_2-1$ twisted sectors with degree shifting numbers
    $\frac{1}{2}+\frac{i}{2d_2}$ for $1\leq i\leq 2d_2-1$. $\{q\}\times \{p'\}$
    gives rise to $2d_1-1$ twisted
    sectors with degree shifting numbers
    $\frac{1}{2}+\frac{i}{2d_1}$ for $1\leq i\leq 2d_1-1$. $\{q\}\times \{q'\}$
    give rise to $4d_1d_2-1$
    twisted sectors with degree shifting numbers $\frac{i}{2d_1}+\frac{j}{2d_2}$ for all $i,j$ except
    $(i,j)=(0,0)$. Using this information, we can write down the ordinary orbifold cohomology
    $$b_{orb}^0=b_{orb}^4=1, b_{orb}^1=b^3_{orb}=6, b^2_{orb}=6$$
    $$b^{\frac{i}{d_1}}_{orb}=b^{\frac{i}{d_2}}_{orb}=1,
    b^{1+\frac{i}{d_1}}_{orb}=b^{1+\frac{i}{d_2}}_{orb}=3,
    b^{2+\frac{i}{d_1}}_{orb}=b^{2+\frac{i}{d_2}}_{orb}=2, \
    1\leq i\leq d_1-1, 1\leq j\leq d_2-1$$
    $$b^{\frac{i}{d_1}+\frac{j}{d_2}}_{orb}=1, 0\leq i\leq 2d_1-1, 0\leq j\leq 2d_2,
    (i,j)\neq (0,0),(d_1,d_2).
    \leqno(3.5.9)$$
    Next, we compute $H^*_{orb, \alpha}$. In this example, the most of the twisted sectors are
    dormant sectors. To find nondormant sectors, recall that $WP(2,2d_1)\times WP(2,2d_2)=
    WP(1,d_1)\times WP(1,d_2)/\Z_2\times \Z_2$. Let $g$ be the
    generator of the first factor and $h$ be the generator of the second factor.
    The fixed points of $g$ are $\{p,q\}\times WP(1,d_2)$. We have two nondormant sectors
    obtained by dividing by the
    remaining action generated by $h$. However, $\gamma(\alpha)_{g,h}=-1$.
    There is no invariant cohomology of $WP(1,d_2)$ under the
    action of $h$ twisted by $\gamma(\alpha)_g$. Hence, these two nondormant twisted sectors
    give no contribution to twisted orbifold cohomology. Their degree shifting numbers are 1.
    For the
    same reason, $WP(1,d_1)\times \{p',q'\}/g$ gives no
    contribution to twisted orbifold cohomology. The fixed point
    locus of $gh$ consists of  4 points which give 4 nondormant sectors. Again, their degree shifting
    numbers are 1. As we saw in the last example, their flat orbifold bundles are nontrivial.
    Hence, they give no contribution to twisted orbifold cohomology.  Therefore,
    the twisted orbifold cohomology is
    $$b_{orb,\alpha}^0=b_{orb,\alpha}^4=1, b_{orb,\alpha}^1=b^3_{orb,\alpha}=2,
    b^2_{orb,\alpha}=2$$
    $$b^{\frac{i}{d_1}}_{orb,\alpha}=b^{\frac{i}{d_2}}_{orb,\alpha}=1,
    b^{1+\frac{i}{d_1}}_{orb,\alpha}=b^{1+\frac{i}{d_2}}_{orb,\alpha}=3,
    b^{2+\frac{i}{d_1}}_{orb}=b^{2+\frac{i}{d_2}}_{orb,\alpha}=2,
    \ 1\leq i\leq d_1-1, 1\leq j\leq d_2-1$$
    $$b^{\frac{i}{d_1}+\frac{j}{d_2}}_{orb,\alpha}=1,\ 0\leq i\leq 2d_1-1, 0\leq j\leq 2d_2,
    (i,j)\neq (0,0),(d_1,d_2).
    \leqno(3.5.10)$$
    \vskip 0.1in
    \noindent
    {\bf Example 3.5.5 $T^6/\Z_4$: } Here, $T^6=\C^3/\bigwedge$, where $\bigwedge$ is the lattice of
    integral points. The generator of $\Z_4$ acts on $T^6$ as
    $$\kappa: (z_1, z_2, z_3)\rightarrow (-z_1, iz_2, iz_3).\leqno(3.5.11)$$
    This example has been studied by D. Joyce \cite{JO}, where he constructed five different
    desingularizations. However, there is no  discrete torsion in the case which induces
    nontrivial orbifold cohomology.

    First of all, the nontwisted sector contributes one generator to
    $H^{0,0}_{orb}, H^{3,3}_{orb}$, 5 generators to $H^{1,1}_{orb},
    H^{2,2}_{orb}$ and 2 generators to $H^{2,1}_{orb},
    H^{1,2}_{orb}$
    The fixed point locus of $\kappa, \kappa^3$ consists of  16 points
    $$\{(z_1, z_2, z_3)+\bigwedge: z_1\in \{0,\frac{1}{2}, \frac{i}{2}, \frac{1}{2}+\frac{i}{2}\},
    z_2, z_3\in \{0, \frac{1}{2}+\frac{i}{2}\}.$$
    These points are fixed by $\Z_4$. Therefore, they generate
    32 twisted sectors in which 16 correspond to the conjugacy
    class $(\kappa)$ and 16 correspond to the conjugacy class
    $(\kappa^3)$. The sector with conjugacy class $(\kappa)$ has
    degree shifting number 1. The sector with conjugacy class
    $(\kappa^3)$ has degree shifting number 2.

    The fixed point locus of $\kappa^2$ is 16 copies of $T^2$,
    given by
    $$\{(z_1, z_2, z_3)+\bigwedge: z_1\in \C, z_2, z_3\in \{0,\frac{1}{2},
    \frac{i}{2}, \frac{1}{2}+\frac{i}{2}\}\}$$
    Twelve of the 16 copies of $T^2$ fixed by $\kappa^2$ are
    identified in pairs by the action of $\kappa$, and these
    contribute 6 copies of $T^2$ to the singular set of $T^6/\Z_4$.
    On the remaining 4 copies $\kappa$ acts as $-1$, so these contribute
    4 copies of $S^2=T^2/\{\pm 1\}$ to the singular set. The degree shifting number
     of these 2-dimensional twisted sectors  is 1.

     Next, we construct local systems. We start with two-dimensional twisted sectors.
     Since $\kappa^{-2}=\kappa^2$, the condition (2) of Definition 2.1 tells us
    that the flat orbifold line bundle $L$ over two-dimensional sectors has the
    property $L^2=1$. Now, we assign the trivial line bundle to all
     $T^2$-sectors and $k (k=0,1,2,3,4)$ $S^2=T^2/\{\pm 1\}$-sectors.
     For the remaining $S^2=T^2/\{\pm 1\}$-sectors, we assign a flat orbifold line
     bundle $T^2\times \C/\{\pm 1\}$. For the zero-dimensional sectors, they are all
    points of two-dimensional sectors. If we assign a trivial bundle on a two-dimensional
    sector, we just assign the trivial bundle to its point sectors. For these
    two-dimensional sectors with nontrivial flat line bundle, we need to be careful to choose
    the flat orbifold line bundle on its point sectors to ensure the condition (3) of
    Definition 2.1. Suppose that $\Sigma$ is one of the 2-dimensional sectors supporting
    a nontrivial flat orbifold line bundle. It contains 4 singular points which generate the
    point sectors. Let $x$ be one of the 4 points. $x$ generates two sectors given by the conjugacy
    classes $(\kappa), (\kappa^3)$. For condition (3),  we have to consider
     the conjugacy class of the triple $(g_1,g_2, g_3)$ with
     $g_1g_2g_3=1$. The only nontrivial choices are
     $(\g)=(\kappa,\kappa, \kappa^2), (\kappa^2, \kappa^3,\kappa^3)$.
    The corresponding components of $X_{(\g)}$ are exactly  these singular points.
      Clearly, $x$ is  fixed by the whole group $\Z_4$. The orbifold line bundle is given
    by the action of $\Z_4$ on $\C$. Consider the component of $X_{(\g)}$ generated by $x$.
    The pull-back of the flat orbifold line bundle from the 2-dimensional sector ($(\kappa^2)$-sector)
    is given by the action
    $\kappa v=-v$. A moment's thought tells us that for sectors $(\kappa), (\kappa^3)$, we should
    assign a flat
    orbifold line bundle given by the action of $\Z_4$ on $\C$ as $\kappa v=iv.$ It is easy to check
    that for the above choices the condition (3) is satisfied for $X_{(\g)}$. Therefore,
    the twisted sectors given by $(x, (\kappa)), (x, (\kappa^3))$ give no contribution
    to twisted orbifold cohomology.  Suppose that the resulting local system is
     $\L_k$. For the sectors with trivial
     line bundle, they contribute $6+k$ generators to $H^{1,1}_{orb},
     H^{2,2}_{orb}$ and 6 generators to $H^{2,1}_{orb},
     H^{1,2}_{orb}$. Its point sectors contribute $4k$ generators to
    $H^{1,1}_{orb}, H^{2,2}_{orb}$.
     The remaining sectors contribute $4-k$ generators to
     $H^{2,1}_{orb}, H^{1,2}_{orb}$. Its point sectors give no contribution.
    Moreover, the nontwisted sector contributes
    $$h^{0,0}=h^{3,3}=2, h^{1,1}=5.$$
    In summary, we obtain
     $$\dim H^{0,0}_{orb}(T^6/\Z_4, \L_k)=\dim H^{3,3}_{orb}(T^6/\Z_4,
     \L_k)=1, \dim H^{1,1}_{orb}(T^6/\Z_4, \L_k)=\dim H^{2,2}_{orb}(T^6/\Z_4,
     \L_k)=11+5k,$$
    $$\dim H^{1,2}_{orb}(T^6/\Z_4, \L_k)= \dim H^{2,1}_{orb}(T^6/\Z_4,
     \L_k)=12-k\leqno(3.5.12)$$
     Our calculation matches the Betti numbers of Joyce's desingularizations \cite{JO}.

\section{Orbifold K-theory}
    It was known classically that any reduced orbifold can be expressed as
    $P/G$, where $P$ is a smooth manifold and $G$ is a compact
    Lie group acting smoothly on $P$ such that $G$ has only finite
    isotropy subgroups. Therefore, it is natural to use the equivariant
    theory of $P$ to capture the theory on $P/G$. The classical
    case is the equivariant cohomology $H^*_G(P)$. If $G$ is
    connected, it is known that $H^*_G(P, \C)=H^*(P/G, \C)$,i.e., the
    nontwisted sector. In the case of a global quotient $X/G$ for a
    finite group $G$ Atiyah-Segal and others showed that
    equivariant $K$-theory $K_G(X)$ carries more information. In
    fact, $K_G(X)\otimes \C= H^*_{orb}(X/G, \C)$. This section has
    two purposes. (1) We would like to generalize Atiyah-Segal \cite{AS} and
    results of others to a general orbifold; (2) more importantly, we
    want to incorporate discrete torsion into our theory. The
    latter leads to some unexpected structure unique from
    the $K$-theory point of view. This section is a joint work with
    Alejandro Adem \cite{AR}. Some related work has been done in
    the context of the $K$-theory of algebraic vector bundles in algebraic geometry
    by Vistoli, B. Toen \cite{T}, and algebraic K-theory of
    $\C^*$-algebra module by Marcolli-Mathai \cite{MM}.

    \subsection{Projective representation}
    Mathematically, our construction is based on projective
    representation. This subsection is a review of basic material
    on projective representations of finite group. Throughout this
    subsection, we will assume that $G$ is finite. Most of the
    background results which we list appear in \cite{KA} Chapter
    III.
    \vskip 0.1in
    \noindent
    {\bf Definition 4.1.1: }{\it Let $V$ denote a finite dimensional complex vector space. A
mapping $\rho: G\to GL(V)$ is called a projective representation
of $G$ if there exists a function $\alpha : G\times G\to \C^*$
such that $\rho (x)\rho (y) = \alpha (x,y)\rho (xy)$ for all
$x,y\in G$ and $\rho (1) = Id_V$.}
    \vskip 0.1in

Note that $\alpha$ defines a $\C^*$--valued cocycle on $G$, i.e.
$\alpha \in Z^2(G,\C^*)$. Also there is a one-to-one
correspondence between projective representations of $G$ as above
and homomorphisms from $G$ to $PGL(V)$. We will be interested in
the notion of \emph{linear equivalence} of projective
representations.
    \vskip 0.1in
    \noindent
    {\bf Definition 4.1.2: }{\it
Two projective representations $\rho_1 : G\to GL(V_1)$ and
$\rho_2: G\to GL(V_2)$ are said to be linearly equivalent if there
exists a vector space isomorphism $f:V_1\to V_2$ such that
$\rho_2(g)=f\rho_1(g)f^{-1}$ for all $g\in G$.}
    \vskip 0.1in

If $\alpha$ is the cocycle attached to $\rho$, we say that $\rho$
is an $\alpha$--representation on the space $V$. We list a couple
of basic results
    \vskip 0.1in
    \noindent
    {\bf Lemma 4.1.3: }{\it
Let $\rho_i$, $i=1,2$ be an $\alpha_i$--representation on the
space $V_i$. If $\rho_1$ is linearly equivalent to $\rho_2$, then
$\alpha_1$ is equal to $\alpha_2$.}
    \vskip 0.1in

It is easy to see that given a fixed cocycle $\alpha$, we can take
the direct sum of any two $\alpha$--representations. Hence we can
introduce
    \vskip 0.1in
    \noindent
    {\bf Definition 4.1.4: }{\it
We define $M_{\alpha}(G)$ as the monoid of linear isomorphism
classes of $\alpha$--representations of $G$. Its associated
Grothendieck group will be denoted $R_{\alpha}(G)$.}
    \vskip 0.1in

In order to use these constructions we need to introduce the
notion of a \emph{twisted group algebra}. If $\alpha :G\times G\to
\C^*$ is a cocycle, we denote by $\C^{\alpha}G$ the vector space
over $\C$ with basis $\{\overline{g}~|g\in G\}$ with product
$$\overline{x}\cdot \overline{y} = \alpha (x,y) \overline{xy}$$
extended distributively.

One can check that $\C^{\alpha}G$ is a $\C$--algebra with
$\overline{1}$ as the identity element. This algebra is called the
\emph{$\alpha$--twisted group algebra} of $G$ over $\C$. Note that
if $\alpha (x,y)=1$ for all $x,y\in G$, then $\C^{\alpha}G=\C G$.
    \vskip 0.1in
    \noindent
    {\bf Definition 4.1.5: }{\it
If $\alpha$ and $\beta$ are cocycles, then $\C^{\alpha} G$ and
$\C^{\beta}G$ are equivalent if there exists a $\C$-algebra
isomorphism $$\psi : \C^{\alpha}G\to \C^{\beta}G$$ and a mapping
$t:G\to \C^*$ such that $\psi (\overline{g}) = t(g) \tilde{g}$ for
all $g\in G$, where $\{\overline{g}\}$ and $\{\tilde{g}\}$ are
bases for the two twisted algebras.}
    \vskip 0.1in

We have a basic result which classifies twisted group algebras.
    \vskip 0.1in
    \noindent
    {\bf Theorem 4.1.6: }{\it
We have an isomorphism between twisted group algebras,
$\C^{\alpha}G\simeq \C^{\beta}G$, if and only if $\alpha$ is
cohomologous to $\beta$; hence if $\alpha$ is a coboundary,
$\C^{\alpha}G\simeq \C G$ Indeed, $\alpha\mapsto \C^{\alpha}G$
induces a bijective correspondence between $H^2(G,\C^*)$ and the
set of equivalence classes of twisted group algebras of $G$ over
$\C$.}
    \vskip 0.1in

Next we recall how these twisted algebras play a role in
determining $R_{\alpha}(G)$.
    \vskip 0.1in
    \noindent
    {\bf Theorem 4.1.7: }{\it
There is a bijective correspondence between
$\alpha$--representations of $G$ and $\C^{\alpha}G$--modules. This
correspondence preserves sums and bijectively maps linearly
equivalent (respectively irreducible, completely reducible)
representations into isomorphic (respectively irreducible,
completely reducible) modules.}

    Recall the Definition 3.1.15 that an element $g\in G$ is said to
be \emph{$\alpha$--regular} if $\alpha (g,x)=\alpha (x,g)$ for all
$x\in C_G(g)$.

Note that the identity element is $\alpha$--regular for all
$\alpha$. Also one can see that $g$ is $\alpha$--regular if and
only if $\overline{g}\cdot
\overline{x}=\overline{x}\cdot\overline{g}$ for all $x\in C_G(g)$.

If an element $g\in G$ is $\alpha$--regular, then any conjugate of
$g$ is also $\alpha$--regular, hence we can speak of
$\alpha$--regular conjugacy classes in $G$. For technical purposes
we also want to introduce the notion of a `standard' cocyle; it
will be a cocycle $\alpha$ with values in $\C^*$ such that (1)
$\alpha (x,x^{-1})=1$ for all $x\in G$ and (2) $\alpha (x,g)\alpha
(xg, x^{-1})=1$ for all $\alpha$--regular $g\in G$ and all $x\in
G$. Expressed otherwise, this simply means that $\alpha$ is
standard if and only if for all $x\in G$ and for all
$\alpha$--regular elements $g\in G$, we have $\overline{x}^{-1} =
\overline{x^{-1}}$ and $\overline{x}\overline{g}\overline{x}^{-1}
= \overline{xgx^{-1}}$. It can be shown that in fact any
cohomology class $c\in H^2(G,\C^*)$ can be represented by a
standard cocycle, hence we will make this assumption from now on.

The next result is basic:
    \vskip 0.1in
    \noindent
    {\bf Theorem 4.1.8: }{\it
If $r_{\alpha}$ is equal to the number of non--isomorphic
irreducible $\C^{\alpha}G$--modules, then this number is equal to
the number of distinct $\alpha$--regular conjugacy classes of $G$.
In particular $R_{\alpha}(G)$ is a free abelian group of rank
equal to $r_{\alpha}$.}
    \vskip 0.1in

    \subsection{Twisted Equivariant K--theory and Decomposition Theorem}
    In this subsection, we assume that $G$ is a semi-direct product of
    a compact Lie group $H$ and a discrete group $\Gamma$
    where $H$ is a compact Lie group and $\Gamma$ is a discrete group.
    Suppose $\alpha \in H^2(\Gamma, U(1))$.  We have a group
    extension

$$1\to S^1\to \tilde{\Gamma}\to \Gamma\to 1.$$
    Let $\tilde{G}$ be the semi-direct product of $H$ and $\tilde{\Gamma}$.

    Suppose that $G$ acts on a smooth manifold $X$ such that
    $X/G$ is compact and the action has only finite isotropy
    subgroup. It is well-known that $Y=X/G$ is an orbifold.
We are now ready to define a twisted version of equivariant
$K$-theory.
    \vskip 0.1in
    \noindent
    {\bf Definition 4.2.3: }{\it
An $\alpha$--twisted $G$--vector bundle on $X$ is a complex vector
bundle $E\to X$ such that $S^1$ acts on the fibers through complex
multiplication so that the action extends to an action of
$\tilde{G}_{\alpha}$ on $E$ which covers the given $G$--action on
$X$.}
    \vskip 0.1in

In fact $E\to X$ is a $\tilde{G}_{\alpha}$--vector bundle, where
the action on the base is via the projection onto $G$ and the
given $G$--action. Note that if we divide out by the action of
$S^1$, we obtain a \emph{projective} bundle over $X$.
    \vskip 0.1in
    \noindent
    {\bf Definition 4.2.4: }{\it
We define the $\alpha$--twisted $G$--equivariant $K$--theory of
$X$, denoted by $^\alpha K_G(X)$, as the Grothendieck group of
isomorphism classes of $\alpha$--twisted $G$--bundles over $X$.}
    \vskip 0.1in

We begin by considering the case $\alpha=0$; this corresponds to
the split extension $G\times S^1$. Any ordinary $G$--vector bundle
can be made into a $G\times S^1$--bundle via scalar multiplication
on the fibrers; conversely a $G\times S^1$--bundle restricts to an
ordinary $G$--bundle. Hence we readily see that $^{\alpha}
K_G(X)=K_G(X)$.

Next we consider the case when $X$ is equal to a point. It is easy
to verify that we obtain $R_{\alpha}(G)$. More generally, if we
consider an orbit $G/H$, then we have $^\alpha
K_G(G/H)=R_{res^G_H(\alpha)}(H)$.

The reader may have noticed that our twisted equivariant
$K$--theory does not have a product structure. Moreover it depends
on a choice of a particular cohomology class in $H^2(\Gamma,S^1)$.
Our next goal is to relate the different twisted versions by using
a product structure inherited from the additive structure of group
extensions.

Suppose we are given $\alpha, \beta$ in $H^2(\Gamma,S^1)$,
represented by central extensions $1\to S^1\to\tilde{\Gamma}_1\to
\Gamma \to 1$ and $1\to S^1\to\tilde{\Gamma}_2\to \Gamma\to 1$.
These give rise to a central extension of the form $$1\to
S^1\times S^1\to \tilde{\Gamma}_1\times\tilde{\Gamma}_2 \to
\Gamma\times \Gamma \to 1.$$ Now we make use of the diagonal
embedding $\Delta:\Gamma\to \Gamma\times \Gamma$ and the product
map $\mu: S^1\times S^1\to S^1$ to obtain a central extension
$$1\to\mu (S^1\times S^1)\to\tilde{\Gamma}\to \Delta (\Gamma)\to
1.$$ This operation corresponds to the sum of cohomology classes,
i.e. the above extension  represents $\alpha + \beta$. Note that
$ker~\mu = \{ (z, z^{-1})\}\subset S^1\times S^1$. Furthermore, we
can pull back the above construction over $G$.

Now consider an $\alpha$--twisted bundle $E\to X$ and a
$\beta$--twisted bundle $F\to X$. Consider the tensor product
bundle $E\otimes F\to X$. Clearly it will have a
$\tilde{G}_1\times\tilde{G}_2$ action on it, which we can restrict
to the inverse image of $\Delta (G)$. Now note that $ker~\mu$ acts
trivially on $E\otimes F$, hence we obtain a $\tilde{G}$ action on
$E\otimes F$, covering the $G$--action on $X$. This is an $\alpha
+ \beta$--twisted bundle over $X$. Hence we have defined a product
$$^\alpha K_G(X)\otimes^\beta K_G(X)\to ^{\alpha +\beta}K_G(X)$$
which prompts us to introduce the following definition.
    \vskip 0.1in
    \noindent
    {\bf Definition 4.2.5: }{\it
The \emph{total twisted equivariant $K$--theory} of a $G$ space
$X$ is defined as

$$TK_G(X) = \bigoplus_{\alpha\in H^2(\Gamma,S^1)}~^\alpha
K_G(X)$$}
    \vskip 0.1in

Using the product above, we deduce that $TK_G(X)$ is a graded
algebra, as well as a module over $K_G(X)$.

We obtain a purely algebraic construction from the above when $X$
is a point. Namely we obtain the \emph{total twisted
representation ring} of $G$, defined as
$$TR(G)=\bigoplus_{\alpha\in H^2(\Gamma, S^1)}~ R_{\alpha}(G),$$
endowed with the graded algebra structure defined above. Note that
if $\alpha$ is the cocycle representing a cohomology class, then
$\alpha^{-1}$ will represent $-\alpha$. Hence we see that
$\rho\mapsto\rho^*$ defines an isomorphism between $R_{\alpha}(G)$
and $R_{-\alpha}(G)$ (indeed, using vector bundles instead we can
easily extend this to show that $^\alpha K_G(X)$ is isomorphic to
$^{-\alpha}K_G(X)$).

     Next, we relate our (twisted) equivariant K-theory to twisted
     orbifold cohomology. Note that $X/H$ is an orbifold covering of
     $X/G$ with covering group $\Gamma$. Therefore, there is a surjective
     homomorphism $\pi^{orb}_1(X)\rightarrow \Gamma$. $\alpha\in H^2(\Gamma, U(1))$
     induces a discrete torsion of $X$ (still denoted by $\alpha$). Then, we have
     following theorem.

    \vskip 0.1in
    \noindent
    {\bf Theorem 4.2.6: }{\it
Suppose that $G$ acts on $X$ such that (i) $X/G$ is compact; (ii)
the action has only finite isotropy subgroups. For any $\alpha\in
H^2(\Gamma,S^1)$ we have a decomposition

$$ ^\alpha K_G(X)\otimes\C \cong H^*_{orb,\alpha}(X/G, \C),$$}
    \vskip 0.1in

    Recall that $H^*_{orb, \alpha}(X/G, \C)$ is a summation over
    each sector. Therefore, Theorem 4.2.6 can be viewed as a
    decomposition theorem of twisted equivariant $K$-theory.
    \vskip 0.1in

    {\bf Proof: }
    We outline a proof in the case that $G$ is a finite group.
    The general case requires a more complicated argument. We refer
    readers to our paper.

    The proof  requires constructing an $\alpha$--twisted equivariant Chern
character.  Fix the class $\alpha\in H^2(G,S^1)$; for any subgroup
$H\subset G$ we let $\alpha_H = res^G_H(\alpha)$. Given any such
$H$, we can associate

$$H\mapsto R_{\alpha_h}(H)$$ Note the special case when $H=<g>$, a
cyclic subgroup. As $H^2(<g>,S^1)=0$, $R_{\alpha_{<g>}}(<g>)$ is
additively isomorphic to $R(<g>)$.

As mentioned before, we can assume that $\alpha$ is a standard
cocycle. If $z\in C_G(g)$, then it will define an action in the
following manner, where we use the $\alpha$--twisted product:

$$\overline{z}\cdot\overline{g}\cdot\overline{z}^{-1} =\alpha
(z,g)\alpha (g,z)^{-1}\overline{g}$$ Recalling our definition of
the character $L_g^{\alpha}$ for $C_G(g)$, we see that it agrees
precisely with $z\mapsto \alpha (z,g)\alpha (g,z)^{-1}$. Hence we
infer from this that there is an isomorphism of $C_G(g)$--modules

$$R_{\alpha_{<g>}}(<g>)\otimes\cong R(<g>) \otimes L_g^{\alpha}$$
for all $g\in G$. Note that $R(<g>)$ is a trivial
$C_G(g)$--module.

Now consider $E\to X$, an $\alpha$--twisted bundle over $X$; it
restricts to an $\alpha_{<g>}$--twisted bundle over the fixed
point set $X_g$. We have isomorphisms of $C_G(g)$--modules:
$$^\alpha K_{<g>}(X_g)\otimes\cong K(X_g)\otimes
R_{\alpha_{<g>}}(<g>)\otimes\cong K(X_g)\otimes R(<g>)\otimes
L_g^{\alpha}.$$ Let $u:R(<g>)\to $ denote the map $\chi\mapsto
\chi (g)$. Then the composition

$$^\alpha K_G(X)\to K(X_g)\otimes R(<g>)\otimes L_g^{\alpha} \to
K(X_g)\otimes L_g^{\alpha}$$ has image lying in the invariants
under the $C_G(g)$--action. Hence we can put these together to
yield a map

$$ ^\alpha K_G(X)\otimes\C \cong \bigoplus_{(g)} (K(X_g)\otimes
L_g^{\alpha})^{C_G(g)}\cong H^*_{orb, \alpha}(X/G, \C)
.\leqno(4.2.1)$$
    One checks
that this induces an isomorphism on orbits $G/H$; the desired
isomorphism follows from using induction on the number of
$G$--cells in $X$ and a Mayer-Vietoris argument (as in \cite{AS}).
$\Box$

\subsection{Orbifold K-theory}
   In last section, we develop the theory in an equivariant setting,
   where our starting point is a smooth manifold with a smooth
   action of a Lie group. Hence, the quotient is an orbifold. However, there are many different way to
   represent an orbifold as such a quotient. In this section, our starting point is the orbifold
   itself. We will use equivariant theory as an intermediate
   object.

   Recall that for every orbifold $X$ there is an orbifold universal
   cover $Y\rightarrow X$ such that the covering group is
   $\pi^1_{orb}(X)$. Then a discrete torsion $\alpha$ is an
   element of $H^2(\pi^1_{orb}(X), S^1)$. Suppose that
   $$1\rightarrow S^1\rightarrow \tilde{G}\rightarrow G\rightarrow
   1$$
   is the extension determined by $\alpha$. We define
   \vskip 0.1in
   \noindent
   {\bf Definition 4.3.1: }{\it We define $^\alpha K_{orb}(X)$ to
   be the Grothendieck group of isomorphism classes of $\alpha$-twisted $\pi^1_{orb}(X)$-orbifold
    bundles over $Y$ and total orbifold K-theory
    $$TK_{orb}(X)=\bigoplus_{\alpha\in H^2(\pi^1_{orb}(X), S^1)} ^\alpha K_{orb}(X).\leqno(4.3.1)$$
    Furthermore, we define the (twisted) orbifold Euler characteristic
    $\chi_{\alpha}(X)=\chi( ^\alpha K_{orb}(X)).$}
    \vskip 0.1in

    Suppose that $X$ is reduced orbifold. So is the orbifold
    universal cover $Y$. Choose a Riemannian metric on $X$. The pull-back metric
    on $Y$ is $\pi^1_{orb}(X)$-invariant. It is well-known that the
    frame bundle $P(Y)$ is a smooth manifold such that $O(n)$ acts
    on $P(Y)$ with finite isotropy subgroups. Since $\pi^1_{orb}(X)$
    acts as isometries, its action lifts to the action on
    $P(Y)$. Moreover, the two actions commute. Then, we take $G=SO(n)\times \pi^{orb}_1(X)$.
    It acts on $P(Y)$ with finite isotropy
    subgroups. It is obvious that $X=P(Y)/G$.  It is easy to check that
    $$^\alpha K_{orb}(X)\cong ^\alpha K_G(P(Y)).\leqno(4.3.2)$$
    Using our decomposition theorem,
    \vskip 0.1in
    \noindent
    {\bf Theorem 4.3.2: }{\it Suppose that $X$ is a reduced
    orbifold. Then, for any discrete torsion $\alpha\in H^2(\pi^1_{orb}(X), S^1)$,
    there is an additive isomorphism
    $$^\alpha K_{orb}(X)\otimes \C\cong H^*_{orb, \alpha}(X,
    \C).\leqno(4.3.2)$$}
    \vskip 0.1in

\subsection{Examples}
    \noindent
    {\bf Example 4.4.1: } Suppose that $X$ is a point with a
    nontrivial group $G$. It is obvious that $^\alpha
    K_G(X)=R_{\alpha}(G)$. Our theorem 4.3.2 yields that
$R_{\alpha}(G)\otimes \C$ has rank equal to the number of
conjugacy classes of elements in $G$ such that the associated
character $L_g^{\alpha}$ is trivial. This of course agrees with
the number of $\alpha$--regular conjugacy classes, as indeed
$^\alpha K_G(X)=R_{\alpha}(G)$. Moreover, the twisted orbifold
Euler charateristic $\chi_{\alpha}$ equals  the number of
$\alpha$-regular conjugacy classes.
    \vskip 0.1in
    \noindent
    {\bf Example 4.4.2: } We will now consider the case of a weighted projective space.
$\C P(d_1,d_2)$ where $(d_1,d_2)=1$. Let $S^1$ act on the unit
sphere in $\C^2$ via

$$(v,w)\mapsto (z^{d_1}v, z^{d_2}w).$$ The space $\C P(d_1,d_2)$
is the quotient under this action, and it has two singular points.
$x=[1,0]$ and $y=[0,1]$. In this case the Lie group used to
present the orbifold is $SO(2)=S^1$ and the corresponding isotropy
subgroups are precisely $\Z/d_1$ and $\Z/d_2$. Their fixed point
sets are disjoint circles in $S^3$, hence the formula for the
orbifold K-theory yields $$K^*_{orb}(\C P(d_1,d_2))\cong_{\C}
\bigoplus_{d_1+d_2-2}K^*(\{ *\})\oplus K^*(\CP(d_1,d_2)).$$ This
is an additive formula which yields the following orbifold Euler
characteristic $$\chi_{orb}(\C P(d_1,d_2))= d_1+d_2. $$
    \vskip 0.1in
    \noindent
    {\bf Example 4.4.3: }
    Let $G(\R)$ denote a semisimple $\Q$--group, and $K$
a maximal compact subgroup. Let $\Gamma\subset G(\Q$) denote an
arithmetic subgroup. Then $\Gamma$ acts on $X= G(\R)/K$, a space
diffeomorphic to euclidean space. Moreover if $H$ is any finite
subgroup of $\Gamma$, then $X^H$ is a totally geodesic
submanifold, hence also diffeomorphic to euclidean space. We can
make use of the Borel-Serre completion $\overline X$. This is a
contractible space with a proper $\Gamma$--action such that the
$\overline{X}^H$ are also contractible (we are indebted to A.Borel
for outlining a proof of this \cite{BO}) but having a compact
orbit space ${\overline X}/\Gamma$. This of course is a basic
geometric restriction on these groups, in particular implying that
an arithmetic group has only finitely many conjugacy classes of
finite subgroups, all of bounded order. Furthermore, using
equivariant triangulation, we can assume that $\overline{X}$ is a
proper, finite  $\Gamma$--CW complex. From this the decomposition
theorem allows us to express the orbifold Euler characteristic of
$X/\Gamma$ in terms of group cohomology:
    \vskip 0.1in
    \noindent
    {\bf Theorem 4.4 }{\it
\[
\chi_{orb}(X/\Gamma) = \sum_{\{ (\gamma)~|~\gamma\in
\Gamma~of~finite~order \}}~~~\chi (Z_{\Gamma}(\gamma))
\]}
\vskip 0.1in

We now illustrate this with a well--known example:

Let $K$ be a totally real number field with ring of integers
$\mathcal O$, and let $\zeta_k$ denote the Dedekind zeta function
associated to $k$. The centralizer of every finite subgroup in
$\Gamma = SL_2(\mathcal O)$ is finite, except for $\pm 1$. Let
$n(\Gamma)$ denote the number of distinct conjugacy classes of
elements of finite order in $\Gamma$. Then applying the above
corollary we obtain $$\chi_{orb}(X/\Gamma) = n(\Gamma) - 2 + 2\chi
(X/\Gamma)$$ Using a formula due to Brown \cite{B} for the regular
Euler characteristic, we obtain the following:
$$\chi_{orb}(X/\Gamma) =n(\Gamma) - 2 +4\zeta_k(-1) + \sum_{(H)}
(2-{\frac{4}{|H|}})$$ where $H$ ranges over all
$\Gamma$--conjugacy classes of maximal finite subgroups in
$\Gamma$.

\section{Orbifold Quantum Cohomology}

    In the last two sections, we discuss the stringy topology of
    orbifold from both the cohomological and $K$-theoretic points of
    view. I hope that I have demonstrated that the stringy topology of
    orbifolds
    is a rich field where geometry, group theory and physics naturally intertwine.
    However, one should view it only as a beginning. For example, there is a whole range of
    theories in differential topology related to transversality theory. These theories are
    very difficult to  generalize to singular manifold. I believe that our stringy
    topology opens a door for such an {\em orbifold differential topology}. $K$-theoretic point
    of view of stringy topology suggests that there must be much more interaction between
    such a classical subject as algebraic topology and string theory. However, it is
    beyond author's expertise.

     Now, we turn to stringy geometry of orbifolds. The stringy geometry is much less
     understood than stringy topology. Since we are primarily
     motivated by physics, a natural question of stringy geometry is whether we can
     quantize the theory. Namely, can we build a theory of
     orbifold quantum cohomology such that orbifold cohomology is
     the classical part of orbifold quantum cohomology. The answer
     is yes. However, we need to restrict
     ourselvs to symplectic or projective orbifolds. It is still an interesting question
     whether
     a projective orbifold is always symplectic.  This section  is
     joint work of the author with W. Chen \cite{CR2}.

     \subsection{Orbifold stable map}

    The most important step of our construction of orbifold
    quantum cohomology is to have an appropriate definition of
    stable map. This is a nontrivial step because a naive
    straightforward generalization is wrong. The key new concept
    is that of good map which we introduced in  section 1. In this
    section, we adapt it to the case where the domain is a nodal Riemann
    surface. Then, an orbifold stable map is defined as an
    ordinary stable map which is good. The central theorem of this
    section is that the moduli space of orbifold stable maps is a
    compact, Housdorff, metrizable space. An equivalent
    formulation of orbifold stable map in algebraic geometry was studied independently
    by D. Abromivich and Vistoli \cite{AV}.

     We first recall
    \vskip 0.1in
\noindent{\bf Definition 5.1.1: }{\it A  nodal Riemann surface
with $k$ marked points is a pair $(\Sigma,\z)$ of a connected
topological space $\Sigma=\bigcup\pi_{\Sigma_\nu}(\Sigma_\nu)$,
where $\Sigma_\nu$ is a smooth complex curve  and
$\pi_\nu:\Sigma_\nu\rightarrow \Sigma$ is a continuous map, and
$\z=(z_1,\cdots,z_k)$ are distinct $k$ points in $\Sigma$ with the
following properties.
\begin{itemize}
\item For each $z\in\Sigma_\nu$, there is a neighborhood of it such that
the restriction of $\pi_\nu:\Sigma_\nu\rightarrow \Sigma$ to this
neighborhood is a homeomorphism onto its image.
\item For each $z\in\Sigma$, we have $\sum_\nu \#\pi_\nu^{-1}(z)\leq 2$.
\item $\sum_\nu \#\pi_\nu^{-1}(z_i)=1$ for each $z_i\in\z$.
\item The number of complex curves $\Sigma_\nu$ is finite.
\item The set of nodal points $\{z|\sum_\nu \#\pi_\nu^{-1}(z)=2\}$ is  finite.
\end{itemize}}\vskip 0.1in

A point $z\in\Sigma_\nu$ is called {\it singular} if
$\sum_\omega\#\pi_\omega^{-1}(\pi_\nu(z))=2$. A point
$z\in\Sigma_\nu$ is said to be a {\it marked point} if
$\pi_\nu(z)=z_i\in\z$. Each $\Sigma_\nu$ is called a {\it
component} of $\Sigma$. Let $k_\nu$ be the number of points on
$\Sigma_\nu$ which are either singular or marked, and $g_\nu$ be
the genus of $\Sigma_\nu$; a nodal curve $(\Sigma,\z)$ is called
{\it stable} if $k_\nu+2g_\nu\geq 3$ holds for each component
$\Sigma_\nu$ of $\Sigma$.

A map $\vartheta:\Sigma\rightarrow\Sigma^\prime$ between two nodal
curves is called as {\it isomorphism} if it is homeomorphism and
if it can be lifted to biholomorphisms
$\vartheta_{\nu\omega}:\Sigma_\nu\rightarrow \Sigma^\prime_\omega$
for each component $\Sigma_\nu$ of $\Sigma$. If $\Sigma$
,$\Sigma^\prime$ have marked points $\z=(z_1,\cdots,z_k)$ and
$\z^\prime=(z_1^\prime,\cdots,z_k^\prime)$ then we require
$\vartheta(z_i) =z_i^\prime$ for each $i$. Let $Aut(\Sigma,\z)$ be
the group of automorphisms of $(\Sigma,\z)$.

Each nodal curve $(\Sigma,\z)$ is canonically associated with a
graph $T_\Sigma$ as follows. The vertices of $T_\Sigma$ correspond
to the components of $\Sigma$ and for each pair of components
intersecting each other in $\Sigma$ there is an edge joining the
corresponding two vertices. For each point $z\in\Sigma$ such that
$\#\pi_\nu^{-1}(z)=2$, there is an edge joining the same vertex
corresponding to $\Sigma_\nu$. For each marked point, there is a
half open edge (tail) attaching to the vertex. The graph
$T_\Sigma$ is connected since $\Sigma$ is connected. We can smooth
out all the nodal points to obtain a smooth surface. Its genus is
called arithmetic genus of $\Sigma$. The arithmetic genus can be
computed by the formula

$$g=\sum_\nu g_\nu +rank H_1(T;\Q).$$

    \vskip 0.1in
\noindent{\bf Definition 5.1.2: }{\it A nodal orbicurve is a nodal
marked  curve $(\Sigma,\z)$ with an orbifold structure on each
component $\Sigma_\nu$ satisfying the following conditions.
\begin{itemize}
\item The singular set $\z_\nu=\Sigma\Sigma_\nu$ (in the sense of orbifolds)
of each component $\Sigma_\nu$ is contained in the set of marked
points and nodal points $z$.
\item If $z_\nu\in\Sigma_\nu$ and $z_\omega\in\Sigma_\omega$ satisfy
$\pi_\nu(z_\nu)=\pi_\omega(z_\omega)$, then the germs of
uniformizing systems at $z_\nu$ and $z_\omega$ are the same. (Here
$\nu$ and $\omega$ may be identical.)
\end{itemize}}\vskip 0.1in

Let a disc neighborhood of $z=\pi_\nu^{-1}(z_i)$ be uniformized by
a branched covering map $z\rightarrow z^{m_i}$, and a disc
neighborhood of $z$ such that $\sum_\omega\pi_\omega^{-1}(\pi_\nu
(z))=2$ be uniformized by $z\rightarrow z^{n_j}$. Here $m_i$ and
$n_j$ are allowed to be equal to one, i.e., the corresponding
orbifold structure is trivial there. We denote the corresponding
nodal orbicurve by $(\Sigma,\z,\m,\n)$, where
$\m=(m_1,\cdots,m_k)$ and $\n=(n_j)$.

An {\it isomorphism} between two nodal orbicurves
$\tilde{\vartheta}: (\Sigma,\z,\m,\n)\rightarrow
(\Sigma^\prime,\z^\prime,\m^\prime,\n^\prime)$ is a collection of
$C^\infty$ isomorphisms $\tilde{\vartheta}_{\nu\omega}$ between
orbicurves $\Sigma_\nu$ and $\Sigma_\omega^\prime$ which induces
an isomorphism $\vartheta:(\Sigma,\z)\rightarrow
(\Sigma^\prime,\z^\prime)$. The {\it group of automorphisms} of a
nodal orbicurve $(\Sigma,\z,\m,\n)$ is denoted by
$Aut(\Sigma,\z,\m,\n)$. It is easily seen that
$Aut(\Sigma,\z,\m,\n)$ is a subgroup of $Aut(\Sigma,\z)$ of finite
index.

\vspace{2mm}

\noindent{\bf Definition 5.1.3: }{\it Let $(X,J)$ be an almost
complex orbifold. An  orbifold stable map is a triple
$(f,(\Sigma,\z),\xi)$ described as follows:
\begin{enumerate}
\item $f$ is a continuous map from a nodal
curve $(\Sigma,\z)$ into $X$ such that each $f_\nu=f\circ\pi_\nu$
is a pseudo-holomorphic map from $\Sigma_\nu$ into $X$.
\item $\xi=\{\xi_\nu\}$ where each $\xi_\nu$ is an orbifold structure on
$\Sigma_\nu$ together with a compactible system of the (unique)
$C^\infty$ lifting of $f_\nu:\Sigma_\nu\rightarrow X$. The
orbifold structure on $\Sigma_\nu$ has its set of orbifold points
contained in $\z_\nu$ where $\z_\nu$ is the union of marked points
and singular points in $\Sigma_\nu$. Moreover, the homomorphism on
the local group is injective.
\item Let $\{\tilde{f}_{\nu,UU^\prime},\lambda_\nu\}$ be the compatible system
defined by $\xi_\nu$ for each $\nu$. Then for any
$z_\nu\in\Sigma_\nu$ and $z_\omega\in\Sigma_\omega$ satisfy
$\pi_\nu(z_\nu)=\pi_\omega(z_\omega)$ (here $\nu$ and $\omega$ may
be identical); if we let $p=f(\pi_\nu(z_\nu))$, and the group
homomorphism of $\lambda_\nu$ at $z_\nu$ be given by $e^{2\pi
i/n_\nu}\rightarrow g_\nu$ and the group homomorphism of
$\lambda_\omega$ at $z_\omega$ be given by $e^{2\pi
i/n_\omega}\rightarrow g_\omega$, then $n_\nu=n_\omega$ and
$g_\nu=g_\omega^{-1}$ in $G_p$.
\item Let $k_\nu$ be the order of the set $\z_\nu$, namely the
number of points on $\Sigma_\nu$ which are singular (i.e. nodal or
marked ); if $f_\nu$ is a constant map, then $k_\nu+ 2g_\nu\geq
3$.
\end{enumerate}
We will call $\xi$ a twisted boundary condition of
$f:(\Sigma,\z)\rightarrow X$. cf. Remark 2.2.10.}\vskip 0.1in

A stable map $(f,(\Sigma,\z),\xi)$ determines  a unique  orbifold
structure $(\Sigma,\z,\m,\n)$ on $(\Sigma,\z)$, as part of
$\xi=\{\xi_\nu\}$. We introduce an equivalence relation amongst
the set of stable maps as follows: two stable maps
$(f,(\Sigma,\z),\xi)$ and
$(f^\prime,(\Sigma^\prime,\z^\prime),\xi^\prime)$ are {\it
equivalent} if there exists an isomorphism $\vartheta:
(\Sigma,\z,\m,\n)\rightarrow
(\Sigma^\prime,\z^\prime,\m^\prime,\n^\prime)$ such that
${f}^\prime\circ\vartheta={f}$, and the compatible systems defined
by $\xi^\prime$ pull back via $\vartheta$ to compatible systems
isomorphic to the ones defined by $\xi$ (we write this as
$\xi^\prime\circ\vartheta=\xi$). The {\it automorphism group} of a
stable map $(f,(\Sigma,\z),\xi)$, denoted by
$Aut(f,(\Sigma,\z),\xi)$, is defined by $$
Aut(f,(\Sigma,\z),\xi)=\{{\vartheta}\in Aut(\Sigma,\z,\m,\n)|
{f}\circ\vartheta={f}, \xi\circ\vartheta=\xi\}. $$

The proof of the following lemma is routine and is left to the
readers.

    \vskip 0.1in
\noindent{\bf Lemma 5.1.4: }{\it The automorphism group of an
orbifold stable map is finite.}\vskip 0.1in

Given a stable map $(f,(\Sigma,\z), \xi)$, there is an associated
homology class $f_\ast([\Sigma])$ in $H_2(X;\Z)$ defined by
$f_\ast([\Sigma])=\sum_\nu(f\circ\pi_\nu)_\ast[\Sigma_\nu]$. On
the other hand, for each marked point $z$ on $\Sigma_\nu$, say
$\pi_\nu(z) =z_i\in\z$, $\xi_\nu$ determines, by the group
homomorphism at $z$, a conjugacy class $(g_i)$ where $g_i\in
G_{f(z_i)}$. We thus have a map $ev$ sending each (equivalence
class of) stable map into $\widetilde{\Sigma_k X}$ by
$(f,(\Sigma,\z),\xi)\rightarrow
((f(z_1),(g_1)),\cdots,(f(z_k),(g_k)))$.   Let $\x=\prod_i
X_{(g_i)}$ be a connected component in $\widetilde{X}^k$.

    \vskip 0.1in
\noindent{\bf Definition 5.2.5: }{\it A stable map
$(f,(\Sigma,\z), \xi)$ is said of {\it type $\x$} if
$ev((f,(\Sigma,\z), \xi))\in \x$. Given a homology class $A\in
H_2(X;\Z)$, we let $\overline{\M}_{g,k}(X,J,A,\x)$ denote the {\it
moduli space of equivalence classes of orbifold stable maps of
genus $g$, with $k$ marked points, and of homology class $A$ and
type $\x$}, i.e., $$
\overline{\M}_{g,k}(X,J,A,\x)=\{[(f,(\Sigma,\z),
\xi)]|g_\Sigma=g,\#\z=k, f_\ast([\Sigma])=A,
ev((f,(\Sigma,\z),\xi))\in\x\}. $$}

The rest of this subsection is devoted to giving a topology on
$\overline{\M}_{g,k}(X,J,A,\x)$ and to proving that the moduli
space is compact when $(X,J)$ is a compact symplectic orbifold or
a projective orbifold.

The set of all isomorphism classes of stable curves of genus $g$
with $k$ marked points, denoted by $\overline{\M}_{g,k}$, is
called the {\it Deligne-Mumford compactification} of the muduli
space $\M_{g,k}$ of Riemann surfaces of genus $g$ with $k$ marked
points (assuming $k+2g\geq 3$). The following differential
geometric description of $\overline{\M}_{g,k}$ is standard.

The moduli space $\overline{\M}_{g,k}$ admits a stratification
which is indexed by the combinatorial types of the stable curves.
More precisely, we can associate a connected graph to each nodal
marked Riemann surface by assigning a vertex with an integer
(genus) to each component, an edge connecting two vertices if the
corresponding components intersect, and a tail to each marked
point.

 Let $g_{\nu}$ be the genus of the component $\nu$ and $k_{\nu}$ be the number of edges and tails
  containing
$\nu$ (we count twice the edges both of whose vertices are $\nu$).
Then the data is required to satisfy $$ k_{\nu}+2g_\nu\geq 3,
\hspace{2mm}\mbox{and} \hspace{2mm} \sum_\nu g_\nu +rank
H_1(T;\Q)=g. $$ Let $Comb(g,k)$ be the set of all such objects
$(T,(g_\nu))$. For each element $(\Sigma,\z)\in
\overline{\M}_{g,k}$, there is an associated element of
$Comb(g,k)$ as follows: we take the graph $T=T_\Sigma$, let
$g_\nu$ be the genus of $\Sigma_\nu$. The set of combinatorial
types $Comb(g,k)$ is known to be of finite order.

There is a partial order $\succ$ on $Comb(g,k)$ defined as
follows. Let $(T,(g_\nu))\in Comb(g,k)$. We consider
$(T_\nu,(g_{\nu\omega})) \in Comb(g_\nu,k_\nu)$ for some of the
vertices $\nu=\nu_1,\cdots,\nu_a$ of $T$. We replace the vertex
$\nu$ of $T$ by the graph $T_\nu$, and join the edge containing
$\nu$ to the vertex $o_\nu(i)$, where $i\in \{1,\cdots,k_\nu\}$ is
the suffix corresponding to this edge. We then obtain a new graph
$\tilde{T}$. The number $\tilde{g}_\nu$ is determined from $g_\nu$
and $g_{\nu\omega}$ in an obvious way.  It is easily seen that
$(\tilde{T},(\tilde{g}_\nu), \tilde{o})$ is in $Comb(g,k)$. We
define  $(T,(g_\nu),o)\succ
(\tilde{T},(\tilde{g}_\nu),\tilde{o})$.

We need some information about the structure of  the
Deligne-Mumford compactification $\overline{\M}_{g,k}$. The
following are well-known

    \vskip 0.1in
\noindent{\bf Fact 5.1.6: }{\it  Let $\M_{g,k}(T,(g_\nu))$ be the
set of all stable curves such that the associated object is
$(T,(g_\nu))$. Then
\begin{itemize}
\item $\overline{\M}_{g,k}$ is a compact complex orbifold which admits a
stratification by finitely many strata; each stratum is of the
form $\M_{g,k}(T,(g_\nu))$.
\item There is a fiber bundle $\U_{g,k}(T,(g_\nu))\rightarrow
\M_{g,k}(T,(g_\nu))$ which has the following property. For each
$x=(\Sigma_x,\z_x)\in \M_{g,k}(T,(g_\nu))$, there is a
neighborhood of $x$ in $\M_{g,k}(T,(g_\nu))$ of the form
$U_x=V_x/G_x$, where $G_x=Aut(\Sigma_x,\z_x)$,  such that the
inverse image of $U_x$ in $\U_{g,k}(T,(g_\nu))$ is diffeomorphic
to $V_x\times\Sigma_x/G_x$. There is a complex structure on each
fiber such that the fiber of $y=(\Sigma_y,\z_y)$ is identified
with $(\Sigma_y,\z_y)$ itself, together with a K\"{a}hler metric
$\mu_y$ which is flat in a neighborhood of the singular points and
varying smoothly in $y$.
\item $\M_{g,k}(T^\prime,(g_\nu^\prime))$ is contained in the
compactification of $\M_{g,k}(T,(g_\nu))$ in $\overline{\M}_{g,k}$
only if $(T,(g_\nu))\succ (T^\prime,(g_\nu^\prime))$.
\item Different strata are patched together in a way which is described in
the following local model of a neighborhood of a stable curve in
$\overline{\M}_{g,k}$. A neighborhood of $x=(\Sigma,\z)$ in
$\overline{\M}_{g,k}$ is parametrized by $$ \frac{V_x\times
B_r(\oplus_z T_{z_\nu}\Sigma_\nu\otimes
T_{z_\omega}\Sigma_{\omega})} {Aut(\Sigma,\z)}, $$ where
$z=\pi_\nu(z_\nu)=\pi_\omega(z_\omega)$ (here it may happen that
$\nu=\omega$) runs over all singular points of $\Sigma$; $B_r(W)$
denotes the ball of radius $r$ of vector space $W$. Each $y\in
V_x$ represents a stable curve $(\Sigma_y,\z_y)$ homeomorphic to
$(\Sigma,\z)$, with a K\"{a}hler metric $\mu_y$ which is flat in a
neighborhood of the singular points. Given $y\in V_x$, for each
element $\varsigma= (\sigma_z)\in\oplus_z
T_{z_\nu}\Sigma_\nu\otimes T_{z_\omega}\Sigma_\omega$ there is an
associated stable curve $(\Sigma_{y,\varsigma},\z_{y,\varsigma})$
obtained as follows. Each component $\Sigma_{\nu}$ of $\Sigma_y$
is given a K\"{a}hler metric $\mu_y$ which is flat in a
neighborhood of singular points. This gives a Hermitian metric on
each $T_{z_\nu}\Sigma_{\nu}$. For each non-zero $\sigma_z\in
T_{z_\nu}\Sigma_{\nu}\otimes T_{z_\omega}\Sigma_{\omega}$, there
is a biholomorphic map $\Psi_{\sigma_z}: T_{z_\nu}\Sigma_{\nu}
\setminus\{0\} \rightarrow
T_{z_\omega}\Sigma_{\omega}\setminus\{0\}$ defined by
$u\otimes\Psi_{\sigma_z}(u)=\sigma_z$. Let $|\sigma_z|=R^{-2}$;
then for sufficiently large $R$, the map
$exp_{z_\omega}^{-1}\circ\Psi_{\sigma_z} \circ exp_{z_\nu}$ is a
biholomorphism between $D_{z_\nu}(R^{-1/2})\setminus
D_{z_\nu}(R^{-3/2})$ and $D_{z_\omega}(R^{-1/2})\setminus
D_{z_\omega}(R^{-3/2})$, where $D_{z_\nu}(R^{-1/2})$ is a disc
neighborhood of $z_\nu$ in $\Sigma_{\nu}$ of radius $(R^{-1/2})$
which is flat assuming $R$ is sufficiently large. We glue
$\Sigma_\nu$ and $\Sigma_\omega$ by this biholomorphism. If
$\sigma_z=0$, we do not make any change. Thus we obtain
$(\Sigma_{y,\varsigma},\z_{y,\varsigma})$. Moreover, there is a
K\"{a}hler metric $\mu_{y,\varsigma}$ on $\Sigma_{y,\varsigma}$
which coincides with the K\"{a}hler metric $\mu_y$  on $\Sigma_y$
outside a neighborhood of the singular points, and varies smoothly
in $\varsigma$. Each $\gamma\in Aut(\Sigma,\z)$ takes
$(\Sigma_y,\z_y)$ to $(\Sigma_{\gamma(y)},\z_{\gamma(y)})$
isometrically, so it acts on $\oplus_z T_{z_\nu}\Sigma_\nu\otimes
T_{z_\omega}\Sigma_\omega$. $\gamma$ induces an  isomorphism
between $(\Sigma_{y,\varsigma},\z_{y,\varsigma})$ and
$(\Sigma_{\gamma(y,\varsigma)},\z_{\gamma(y,\varsigma)})$, which
is also an isometry.
\end{itemize}}\vskip 0.1in

Now we define a topology on the moduli space
$\overline{\M}_{g,k}(X,J,A,\x)$. We put a Hermitian metric $h$ on
$(X,J)$ and the  distance function on $X$ is assumed to be induced
from $h$.

\vspace{2mm}

\noindent{\bf Definition 5.1.7: }{\it A sequence of equivalence
classes of stable maps $x_n$ in $\overline{\M}_{g,k}(X,J,A,\x)$ is
said to converge to $x_0\in\overline{\M}_{g,k}(X,J,A,E)$ if there
are representatives $(f_n,(\Sigma_n,\z_n),\xi_n)$ of $x_n$ and a
representative $(f_0,(\Sigma_0,\z_0),\xi_0)$ of $x_0$  the
following conditions hold.
\begin{itemize}
\item For each $n$ (including $n=0$),
there is a set of distinct {\it regular} points
$\{z_{n,1},\cdots,z_{n,a}\}$ (it may happen that this set is
empty) on $\Sigma_n$ which is disjoint from the marked point set
$\z_n$ such that after adding this set to $\z_n$ we obtain a
stable curve in $\overline{\M}_{g,k+a}$, denoted by
$(\Sigma_n,\z_n)^+$. Let $(f_n^+,(\Sigma_n,\z_n)^+,\xi_n^+)$ be
the sequence of stable maps naturally obtained.
\item The sequence $(\Sigma_n,\z_n)^+$ converges to
$(\Sigma_0,\z_0)^+$ in $\overline{\M}_{g,k+a}$. This means that
for sufficiently large $n$, $(\Sigma_n,\z_n)^+$ is identified with
$(\Sigma_{y_n,\varsigma_n},\z_{y_n,\varsigma_n})$ for some
$(y_n,\varsigma_n)$ in the canonical model of a neighborhood of
$(\Sigma_0,\z_0)^+$. Let $\varsigma_n$ be given by
$(\sigma_{z,n})$ and $|\sigma_{z,n}|=R_{z,n}^{-2}$ (here $R_{z,n}$
is allowed to be $\infty$), For each $\mu>\max_z(R_{z,n}^{-1})$ we
put $$ W_{z,n}(\mu)=(D_{z_\nu}(\mu)\setminus
D_{z_\nu}(R_{z,n}^{-1}))\cup (D_{z_\omega}(\mu)\setminus
D_{z_\omega}(R_{z,n}^{-1})),\hspace{2mm}\mbox{and} \hspace{2mm}
W_n(\mu)=\cup_z W_{z,n}(\mu). $$ Then the following holds. First,
for each $\mu>0$, when $n$ is sufficiently large, the restriction
of $\tilde{f}^+_n$ to $\Sigma_{y_n,\varsigma_n}\setminus W_n(\mu)$
converges to $\tilde{f}^+_0$ in the $C^\infty$ topology as a
$C^\infty$ map with an isomorphism class of compatible systems.
Secondly, $\lim_{\mu\rightarrow 0}\limsup_{n\rightarrow\infty}
Diam(f_n(W_{z,n}(\mu)))=0$ for each singular point $z$ of
$\Sigma_0$.
\end{itemize}}\vskip 0.1in

 \noindent{\bf Proposition 5.1.8: }{\it Suppose $X$ is either
a symplectic orbifold with a symplectic form $\omega$ and an
$\omega$-compatible almost complex structure $J$, or a projective
orbifold with an integrable almost complex structure $J$. Then the
moduli space $\overline{\M}_{g,k}(X,J,A,\x)$ is compact and
metrizable.} \vskip 0.1in

 For the proof, the reader is referred to \cite{CR2}.

 \subsection{Orbifold Gromov-Witten Invariants}

For any component $\x=(X_{(g_1)},\cdots,X_{(g_k)})$, there are $k$
evaluation maps.
 $$ e_i:\overline{\M}_{g,k}(X,J,A,\x)\rightarrow
X_{(g_i)}, \hspace{4mm} i=1,\cdots, k. \leqno(5.1) $$ For any set
of cohomology classes $\alpha_i\in
H^{*-2\iota_{(g_i)}}(X_{(g_i)};\Q)\subset H^*_{orb}(X;\Q)$,
$i=1,\cdots,k$, the orbifold Gromov-Witten invariant is defined as
the virtual integral $$ \Psi^{X,J}_{(g,k,A,\x)}(\alpha^{l_1}_1,
\cdots, \alpha^{l_k}_k)=\int^{vir}_{
\overline{\M}_{g,k}(X,J,A,\x)}\prod_{i=1}^k c_1(L_i)^{l_i}e^*_i
\alpha_i,\leqno(5.2) $$ where $L_i$ is the line bundle defined by
the cotangent space of the $i$-th marked point.

\vspace{2mm}

When $g=0$ and $A=0$, the moduli space
$\overline{\M}_{g,k}(X,J,A,\x)$ admits a very nice and elementary
description, based on which we gave an elementary construction of
genus zero, degree zero orbifold Gromov-Witten invariants in
\cite{CR1}. Even in this case, virtual integration is needed where
there is an obstruction bundle. The orbifold cup product (cf.
Theorem 2.3) is defined through these orbifold Gromov-Witten
invariants. In the general case, we need to use the full scope of
the virtual integration machinary developed by \cite{FO},
\cite{LT}, \cite{R2} and \cite{Sie}.

\vspace{2mm}

Singularities of an orbifold impose additional difficulties in
carrying  out virtual integration  in the orbifold case. Due to
the presence of singularities, even on a closed orbifold the
function of injective radius of the exponential map does not have
a positive lower bound. As a consequence, it is not known that a
neighborhood of a (good) $C^\infty$ map into an orbifold can be
completely described by $C^\infty$ sections of the pull-back
tangent bundle via the exponential map. Our approach is a
combination of techniques developed in the smooth case and some
additional techniques for orbifolds. \vspace{2mm}

 The main results of this work are summarized in the following

\vspace{2mm}

\noindent{\bf Theorem 5.2.1:}\hspace{2mm}{\it Let $X$ be a closed
symplectic or projective orbifold. The orbifold Gromov-Witten
invariants defined in (5.2) satisfy the quantum cohomology axioms
of Witten-Ruan for ordinary Gromov-Witten invariants (cf.
\cite{R3}) except that in the Divisor Axiom, the divisor class is
required to be in the nontwisted sector (i.e. in $H^2(X;\Q)$). In
the formulation of axioms, the ordinary cup product is replaced by
the orbifold cup product $\cup_{orb}$ (cf. Theorem 2.3). }

\vspace{2mm}

As a consequence, we have

\vspace{2mm}

\noindent{\bf Theorem 5.2.2:}\hspace{2mm}{\it Let $X$ be a closed
symplectic or projective orbifold. With suitable coefficient ring
$\cal{C}$, the small quantum product and the big quantum product
are well-defined on $H^\ast_{orb}(X;\Q)\otimes\cal{C}$, and have
properties similar to those of the ordinary quantum cohomology. }

\subsection{Examples}
    It is generally a difficult problem to compute orbifold quantum
    cohomology. Much machinery has been developed to compute ordinary
    quantum cohomology. The most important ones are localization
    and surgery techniques. They should have their counterparts in
    orbifold quantum cohomology as well. The problem is that this
    subject is so young that there has not been enough time to develop
    all the machinery. Here, we compute some simple examples by
    direct computations.
    \vskip 0.1in
    \noindent
    {\bf Example 5.3.1: } Let's consider weighted projective space
    $WP(1,1,2,2,2)$. This is an important example in mirror
    symmetry. It has a twisted sector $WP(0,0,2,2,2)$ with local
    group $\Z_2$. Let $\tau$ be the generator. Its degree shifting
    number is one. Recall that $WP(0,0,2,2,2)$ is isomorphic to
    $\P^2$. Therefore,
    $$H^{2i}_{orb}(WP(1,1,2,2,2), \C)=
    \left\{\begin{array}{ll}
    \C & i=0\\
    \C\oplus \C& 1\leq i\leq 3\\
    \C & i=4
    \end{array} \right.$$
    All others are zero.
     Let
    $D_i$ be the hyperplane divisor where $i$-th homogeneous coordinate is
    zero. The first Chern class  $C_1=D_1+D_2+D_3+D_4+D_5$ with
    $2D_1=2D_2=D_3=D_4=D_5$. The weighted projective space is much more
    complicated than projective space. For example, there are three
    types
    of lines in $WP(1,1,2,2,2)$, depending on whether it is in the twisted sector,
    intersects the
    twisted sector transversely or is disjoint from the twisted sector. Their examples
    are $\{[0,0,t,s,0]\}, \{[0,t,s,0,0]\}, \{[t,s,0,0,0]\}$. Let $A_i$ be  its fundamental class.
    It is easy to calculate $C_1(A_1)=2, C_1(A_2)=4, C_1(A_3)=8$. Let's compute the
    invariant for $A_2$; the other two are more difficult to compute. The second kind of line has an
    orbifold point. Consider the moduli space of orbifold stable spheres with
    orbifold points of order $(1,2)$. The complex dimension of the moduli space is
    $6$. Now, we choose a point class $\alpha$ from nontwisted sector and a
    point class $\beta$ from the twisted sector. We are interested in computing the orbifold GW-invariant
    $\Psi_{(A_2, 0)}(\alpha, \beta)$. We observe that there are only two kinds of holomorphic curves
    with homology class $A_2$, a line intersecting twisted sector transversely or a conic in twisted
    sector. A conic in the twisted sector does not pass through a point not in the twisted sector. Therefore,
    for the purpose of calculating the invariant, we only have to consider the lines intersecting transversely
    with the twisted sector. We choose our two points as $[0,1,0,0,0], [0,0,1,0,0]$. There is only one
    such  line passing through these two points. Therefore, there is only one orbifold holomorphic map covering
    the line. This one is regular and hence good. Therefore, $\Psi_{(A_2, 0)}(\alpha,
    \beta)=1$.

\section{Orbifold String Theory Conjectures}
    The physicists believe that orbifold string theory is
    equivalent to ordinary string theory of its
    desingularizations. This belief motivated a body of
    conjectures which we call the {\em Orbifold string theory
    conjecture}. At  present, the case without discrete torsion
    is much better understood than the case with discrete torsion.
    Classical theory is better understood than quantum theory.
    Therefore, we shall start from the classical theory without
    discrete torsion and then discuss the case with discrete
    torsion. We finish  by discussing the case of quantum
    theory. In the physics literature, physicists concern only with
    3-dimensional Calabi-Yau orbifolds. But it is clear that
    much more is true beyond 3-dimensional Calabi-Yau orbifolds.
    Here, we restrict ourselvs to the case of Gorenstein
     reduced orbifolds $X$. In this case, all the local groups
    are subgroup of $SL_2(\C)$ (i.e., a $SL$-orbifold) and hence have integral degree
    shifting number.
    \vskip 0.1in
    \noindent
\subsection{ Classical Case Without Discrete Torsion}
    \vskip 0.1in
     Recall
    \vskip 0.1in
    \noindent
    {\bf Definition 6.1: }{\it A  deformation of $X$ is a
    triple $\pi: U\rightarrow \Delta$ such that $\Delta$ is a
    disc around the origin, $\pi$ is holomorphic and $X=\pi^{-1}(0)$.
    We call $(\pi, U, \Delta)$ is a K\"{a}hler or projective deformation
    if $U$ is K\"{a}hler or projective. We call $X$ {\em the central fiber}
    and $X_t=\pi^{-1}(t)$ a generic fiber. Suppose that $Z$ is a singular space
    such that $K_Z$ is Cartier. A crepant resolution $Y$ of $Z$ is
    a morphism $p:Y\rightarrow Z$ such that $Y$ is smooth and
    $p^*K_Z=K_Y$. A desingularization $Y$ of $X$ is defined as a
    crepant resolution of a generic fiber $X_t$ of a deformation
    of $X$.}
    \vskip 0.1in
    Two extreme cases are $(i)$ $Y$ is a crepant resolution of
    $X$; (ii) The generic fiber $X_t$ is smooth. In this case, we
    call $Y=X_t$ a smoothing of $X$.
    \vskip 0.1in
    \noindent
    {\bf K-Orbifold String Theory Conjecture: }{\it Suppose that
    $X$ is a Gorenstein orbifold and $\pi: Y\rightarrow X$ is a
    crepant resolution. Then there is a natural isomorphism
    between $K_{orb}(X)$ and $K(Y)$.}
    \vskip 0.1in
    Many weaker forms of this conjecture have been studied
    intensively in literature under the name of the McKay
    correspondence. For example, we can replace $K$-theory by the
    Euler number, which we call the $E$-Orbifold string theory
    conjecture. One can also consider a weaker version of
    $K$-orbifold string theory conjecture by dropping naturality.
    Namely, we only consider the corresponding dimensions. We
    label it as $WK$-orbifold string theory conjecture. The best
    result in this direction so far is Batyrev's proof of the
    $WK$-orbifold string theory conjecture for global quotients \cite{B2}.
    Batyrev used a number theoretic method called motivic
    integration invented by Kontsevich \cite{KO}. Actually, Batyrev proved a
    stronger result of the equivalence of Hodge numbers. However, this
    method does not yield a natural map. Moreover, the conjecture
    is completely open for general orbifolds. The hardest part of
    this conjecture is to get a natural map between $K_{orb}(X)$
    and $K(Y)$. Such a map is necessary for us to compare orbifold quantum
    cohomology to quantum cohomology of $Y$. One of the difficulties in constructing such a map is
    that the projection $Y\rightarrow X$ does not pull back class
    from twisted sector. I believe that another formulation of
    orbifold cohomology is needed here and topological methods may
    play an important role. If we go beyond the Gorenstein
    orbifold, orbifold cohomology is rationally graded.
    Batyrev defined  string theoretic Hodge numbers in terms of
    its resolution (not necessarily crepant). The generating
    function of his string theoretic Hodge numbers is not
    necessarily a polynomial, which echoed the rationality of
    grading of orbifold cohomology. It would be an interesting
    question to investigate their relation. In the meantime, we have very
    few examples of non-global quotients which we have calculated.
    It is also very important to calculate more examples.
    Calabi-Yau hypersurface of simplicial toric varieties will be
    a good place to start. The case of complex dimension three has
    been calculated recently by M. Poddar \cite{P} which gives further evidence
    to the $K$-orbifold string theory conjecture.

    The example attracting a lot of attention is the symmetric product
    of the algebraic surface, which is also the best understood example.

    \vskip 0.1in
    \noindent
    {\bf Example 6.2: } Let $S$ be an algebraic surface and $X_k=Sym^k
    (S)$ be the $k$-fold symmetric product of $S$. Then, $X_k$ is a
    $SL$-orbifold and $K_{X_k}$ is Cartier. It is well-known that
    the Hilbert scheme $Hilb_k(S)$ of points of length $k$ is a crepant
    resolution of $X_k$. This case has been extensively studied
    and occupies a special place in  stringy geometry and
    topology. Many years ago, G\"{o}ttsche calculated the generating
    function of the Euler number of $Hilb_k(S)$ and discovered a
    modularity property. Motivated by orbifold string theory,
    Vafa-Witten \cite{VW} calculated the orbifold cohomology group of
    $X_k$ and showed that $\H=\oplus_k H^*_{orb}(X_k, \C)$ is a "Fock
    space" of orbifold conformal field theory. This means that
    $\H$ is an representation of Heisenberg algebra. Under the
    framework of conformal field theory, the Euler characteristic is
    the genus one correlation function and hence modular by
    definition. The $K$-orbifold string theory conjecture predicts
    that $\oplus_k Hilb_k(S)$ is also a representation of
    Heisenberg
    algebra. This was verified by Nakajima \cite{N}.
    Currently, the ring structure of $Hilb_k(S)$ is still unknown
    and is a hot topic right now. Recently, Lehn-Sorger calculated
    ring structure of $Hilb_k(\C^2)$ \cite{LS}. Combined with
    example 3.5.2, it implies that $Hilb_k(\C^2)$ and $Sym_k(\C^2)$ have
    the same ring structure. In the general case, we conjecture
    \vskip 0.1in
    \noindent
    {\bf Conjecture 6.3: }{\it Suppose that $Hilb_k(S)$ has
    hyperk\"{a}hler structure. Then $Hilb_k(S)$ and $Sym_k(S)$ have
    the same ring structure. In general, if $X$  and its crepant resolution $Y$
    have
    hyperK\"{a}hler structure, $Y$ and $X$ have the same ring structure.}
    \vskip 0.1in
    Examples includes $T^4, K3$ and many other open manifolds such as $\C^2$.
    This conjecture is a consequence of the quantum version of
    the orbifold string theory conjecture.

    \vskip 0.1in
    \noindent
\subsection{Classical Case With Discrete Torsion}
    The geometry of orbifold cohomology with discrete torsion is
    much less understood. Originally, Vafa-Witten proposed to use
    it to identify the cohomology of a general desingularization
    which may not be a straightforward crepant resolution. D.
    Joyce \cite{JO} showed that this proposal fails badly to count
    the desingularizations of some simple examples such as
    $T^6/\Z_4$, where there is no discrete torsion but there are
    at least five different desingularizations. To overcome this
    difficulty, I introduced a more general notion of {\em inner
    local system} and constructed a twisted orbifold cohomology
    for any inner local system. Inner local systems were
    successfully used to count all the examples Joyce constructed.
    We propose the following conjecture:
    \vskip 0.1in
    \noindent
    {\bf D-Orbifold String Theory Conjecture: }{\it Suppose that
    $Y$ is a desingularization of $X$. Then there is an inner
    local system $\L$ such that the ordinary cohomology group of $Y$
    is a direct sum of $H^*_{orb}(X, \L)$ and the cohomologies
    generated by the exceptional locus corresponding to small resolution. In particular,
    if $Y$ is a smoothing, its cohomology is the same as twisted orbifold
    cohomology for an inner local system.}
    \vskip 0.1in
    \noindent
    {\bf Remark 6.4: }{\it A generic fiber $X_t$ may not be smooth in
    general. To obtain a crepant resolution, a small resolution
    may be needed. Here a small resolution is a resolution where
    the
    exceptional locus is of codimension at least two.}
    \vskip 0.1in
    The best known example of small resolution is the small resolution
    of a nodal point in three dimension, where the exceptional locus is a
    rational curve of the normal bundle ${\bf O}(-1)+{\bf O}(-1)$.

    \vskip 0.1in
    \noindent
    {\bf Remark 6.5: }{\it The inverse of the D-orbifold string theory
    conjecture is obviously false. We can easily construct
    an example with nontrivial discrete torsion, where the singularities
    are
    of codimension $\geq 3$. In this case, there is no
    deformation.}
    \vskip 0.1in
    It is also an unknown question to identify  desingularizations
    corresponding to discrete torsion.

    It is also very important to understand more examples. Recall
    that symmetric group has a nontrivial discrete torsion. It
    would be an interesting problem to find out if there is a
    desingularization of $Sym_k(S)$ realizing the twisted orbifold
    cohomology.
    \vskip 0.1in
    \noindent
\subsection{Quantum Case}
    \vskip 0.1in
    The ultimate goal of the orbifold string theory conjecture is to
    compare orbifold quantum cohomology of $X$ with the quantum
    cohomology of $Y$. I do not know a precise statement to which I
    could not find a counterexample. However, it is a very useful
    general goal to motivate other better formulated conjectures.
    I do not know how to twist orbifold quantum cohomology using
    an inner local system or a discrete torsion. Therefore, we focus on
    the case without discrete torsion.
    \vskip 0.1in
    \noindent
    {\bf Q-Orbifold String Theory conjecture: }{\it Suppose that
    $Y$ is a crepant resolution of $X$ and $\pi:
    K_{orb}(X)\rightarrow K(Y)$ is the natural additive
    isomorphism given by the $K$-orbifold string theory conjecture.
    Then $\pi$ induces an isomorphism of orbifold quantum
    cohomology up to a mirror transformation.}
    \vskip 0.1in
    The mirror transformation here is similar to the one appearing
    in mirror symmetry \cite{CK}. An interesting case is when $Y$
    is a hyperk\"{a}hler manifold. In this case, there are no quantum corrections
    and the quantum cohomology of     $Y$ is the same as ordinary cohomology.
    Our $Q$-orbifold string theory conjecture becomes a statement
    for the orbifold cohomology of $X$. This is the origin of
    conjecture 6.3.

    The physical prediction in the quantum case is very imprecise. Basically, there
    is a family of superconformal field theories containing both orbifold theory and
    the theory of the resolution. Physics predicts that there is an analytic continuation
    from one theory to  the other. This tells very little about the precise relation between
    them. I believe that orbifold quantum cohomology is different from quantum cohomology
    of crepant resolutions  in general and a mirror transformation is needed.  Actually,
    orbifold cohomology should be naturally related
    to relative quantum cohomology. Suppose that $Z$ is the exceptional divisor of
    the projection $Y\rightarrow X$. We want to
    identify orbifold GW-invariants of $X$ with the relative
    GW-invariants of the pair $(Y, Z)$ introduced by Li-Ruan
    \cite{LR}. Then we can relate relative GW-invariants of $(Y,
    Z)$ with ordinary GW-invariants of $X$. In fact, a
    generalization of Li-Ruan's surgery technique to the orbifold
    category should be very useful for this purpose.

\subsection{Generalization of Orbifold String Theory Conjecture}

   Note that resolution is a special class of birational maps. It is
   natural to recast orbifold string theory conjectures in the context of
   birational geometry. Several years ago, Batyrev \cite{B3} and Wang
   \cite{W} proved that smooth
   birational minimal models have the same Betti numbers. In fact,
   their results are slightly more general.
   \vskip 0.1in
   \noindent
   {\bf Definition 6.6: }{\it $X, Y$ are called $K$-equivalent if
   there is a common resolution $\phi: Z\rightarrow X, \psi:
   Z\rightarrow Y$ such that $\phi^* K_X=\psi^* K_Y$.}
   \vskip 0.1in
   Batyrev and Wang proved that $K$-equivalent smooth
   projective manifolds have the same Betti numbers. At the
   same time, An-Min Li and I \cite{LR} proved that a smooth flop in three
   dimensions
   induces an isomorphism of quantum cohomology. These two results
   inspired the author to propose \cite{R1}
   \vskip 0.1in
   \noindent
   {\bf Quantum minimal model conjecture: }{\it Smooth birational
   minimal models have isomorphic quantum cohomology.}
   \vskip 0.1in

   Here, we observe that if $\phi: Y\rightarrow X$ is a crepant
   resolution, then $Y, X$ are $K$-equivalent. Therefore, we can
   combine orbifold string theory conjectures with quantum minimal
   model conjecture to formulate even more general conjectures.
   \vskip 0.1in
   \noindent
   {\bf $K,Q$-Conjectures: }{\it Suppose that $Y, X$ are
   $K$-equivalent orbifolds. The same statements for
   $K,Q$-orbifold string theory conjectures are true.}
   \vskip 0.1in
   I feel that in the more general context of birational geometry  our $K, Q$-conjectures give us a better
   understanding of relations between orbifold string theory and
   birational geometry than the orbifold string theory conjectures
   itself.
\subsection{Orbifold Mirror Symmetry and Mirror Symmetry in Higher
Dimension}

    One of my first impressions of mirror symmetry is how lucky we
    are in three dimensions where every Calabi-Yau orbifold has a
    crepant resolution. Hence, we have the luxury to consider
    smooth Calabi-Yau 3-folds only. In higher dimensions, it is no
    longer true that every Calabi-Yau orbifold has a crepant
    resolution. We are stuck with orbifolds. Therefore, it makes
    sense to talk about orbifold mirror symmetry,i.e., mirror
    symmetry among orbifolds. Moreover, if one looks at the physical
    literature, physicists clearly consider orbifolds as well.
    For example, the first physical proof of mirror
   symmetry by Greene and Plesser \cite{GP} was based on some
   orbifold model. It is clear that orbifold theory is central to
   mirror symmetry. However, during the mathematization of mirror
   symmetry, orbifold model was replaced by its crepant
   resolution. I tried some simple examples such as Calabi-Yau
   hypersurfaces of weighted projective space $WP(1,1,2,2,2)$. It is
   not clear how mirror symmetry predicts its orbifold quantum
   cohomology. I believe that it is important  to do
   some soul-searching on the role of orbifold string theory in
   mirror symmetry. This is related to another important question
   of
   generalizing  mirror symmetry to higher dimensions, where crepant
   resolution  no longer exists. However, it  still makes
   perfect sense to talk about mirror symmetry between Calabi-Yau
   orbifolds.

\subsection{Other Problems}
   Orbifold String Theory Conjectures or our general $K,
   Q$-conjectures are certainly outstanding problems in stringy geometry and
   topology. There are other very interesting problems in this
   subject as well. Here are some of my favorite problems.

   \vskip 0.1in
   \noindent
   {\bf (1). Orbifold Differential Topology: } As I remarked in
   (3.4.11), it is a very interesting problem to establish a homology
   theory reflecting orbifold  cohomology. This should have a
   profound impact on transversality theory of singular spaces. It
   should also be useful for stack theory in algebraic geometry.

   \vskip 0.1in
   \noindent
   {\bf (2). Relation between two product structures: }
   We have two different products on orbifold cohomology
   groups coming from cohomology and $K$-theory. Each seems to
   reflect one aspect of orbifold theory. The relation between is still mysterious
   to me. It is certainly worth further
   investigation.
   \vskip 0.1in

    Finally, we remark that we only talked about the part of stringy
    geometry and topology motivated by so called type IIA, IIB
    orbifold string theory. There are other types of physical
    orbifold theories which also generate interesting mathematics.
    Unfortunately, the author has very little understanding of
    other types of physical orbifold  theories. I apologize for the
    omission. However, I would like to mention the orbifold
    Landau-Ginzburg theory \cite{KU} and the orbifold elliptic genus
    \cite{BL}
    (orbifold heterotic string theory). Obviously, there is much
    more rich mathematics waiting for us to explore!

    \end{document}